\documentclass[11pt,reqno]{amsart}
\usepackage{amsmath}
  \usepackage{paralist}
  \usepackage{graphics} 
  \usepackage{epsfig} 
\usepackage{graphicx}  
\usepackage{epstopdf}
\usepackage{amssymb}
\usepackage{comment}

 \usepackage[colorlinks=true]{hyperref}
\hypersetup{urlcolor=blue, citecolor=red}

\newtheorem{theorem}{Theorem}[section]
\newtheorem{corollary}{Corollary}

\newtheorem{proposition}{Proposition}

\theoremstyle{definition}

\newtheorem{remark}{Remark}

\newcommand\restr[2]{{
		\left.\kern-\nulldelimiterspace 
		#1 
		\vphantom{\big|} 
		\right|_{#2} 
}}

\title[Inverse problems for the perturbed polyharmonic operator] 
      {Determining rough first order perturbations of the polyharmonic operator}
      
 
\begin{document}
	
\author[Yernat M. Assylbekov]{Yernat M. Assylbekov}
\address{Department of Computational and Applied Mathematics, Rice University, Houston, TX 77005, USA}
\email{yernat.assylbekov@gmail.com}

\author[Karthik Iyer]{Karthik Iyer}
\address{The Vanguard Group, Malvern, PA 19335, USA} \email{karthik2@gmx.com}
	
	\maketitle
	
	\begin{abstract}
		We show that the knowledge of Dirichlet to Neumann map for rough $A$ and $q$ in $(-\Delta)^m +A\cdot D +q$ for $m \geq 2$ for a bounded domain in $\mathbb{R}^n$, $n \geq 3$ determines $A$ and $q$ uniquely. This unique identifiability is proved via construction of complex geometrical optics solutions with sufficient decay of remainder terms, by using property of products of functions in Sobolev spaces.
	\end{abstract}

	\section{Introduction and Statement of Results}
	\subsection{Introduction}
	Let $\Omega \subset \mathbb{R}^n$, $n \geq 3$ be a bounded open set with $C^{\infty}$ boundary. Consider the polyharmonic operator $(-\Delta)^m$ where $m \geq 1$ is an integer. The operator $(-\Delta)^m$ is positive and self-adjoint on $L^2(\Omega)$ with domain $H^{2m}(\Omega) \cap H^m_{0}(\Omega)$, where
	$H^m_{0}(\Omega) = \{ u \in H^m(\Omega): \gamma u = 0 \}$.
	
	This operator can be obtained as the Friedrichs extension starting from the space of test functions; see, for example, \cite{Grubb}. Here and in what follows, $\gamma$ is the Dirichlet trace operator
	$$
	\gamma: H^m(\Omega) \to \prod_{j=0}^{m-1} H^{m-j-1/2}(\partial \Omega),\quad \gamma u = (\restr{u}{\partial \Omega}, \restr{\partial_{\nu}u}{\partial \Omega}, ...., \restr{\partial^{m-1}_{\nu} u}{\partial \Omega}),
	$$
	where $\nu$ is the unit outer normal to the boundary $\partial \Omega$, and $H^s(\Omega)$ and $H^s(\partial \Omega)$ are the standard $L^2$ based Sobolev spaces on $\Omega$ and $\partial \Omega$ respectively for $s \in \mathbb{R}$.
	
	Let us first consider the perturbed polyharmonic operator $\mathcal{L}_{A,q} = (-\Delta)^m + A\cdot D +q$ where $A$ and $q$ are sufficiently smooth and $D = -i \nabla$. For $f=(f_0,f_1,...,f_{m-1}) \in \prod_{j=0}^{m-1} H^{2m - j - 1/2}(\partial \Omega)$, consider the Dirichlet problem
	\begin{equation}\label{1.1}
	\mathcal{L}_{A,q} u = 0\text{ in }\Omega\quad\text{and}\quad \gamma u = f \text{ on }\partial \Omega.
	\end{equation}
	If $0$ is not in the spectrum of $\mathcal{L}_{A,q}$ it can be shown that the Dirichlet problem (\ref{1.1}) has a unique solution $u \in H^{2m}(\Omega)$. We can then define the Dirichlet-to-Neumann $\mathcal{N}_{A,q}$ map  as 
	\begin{align*}
	\mathcal{N}_{A,q}f = (\restr{\partial_{\nu}^m u}{\partial \Omega},...,\restr{\partial_{\nu}^{2m-1} u}{\partial \Omega}) = \tilde \gamma u \in \prod_{j=m}^{2m-1} H^{2m - j - 1/2}(\partial \Omega).
	\end{align*}
	The inverse boundary value problem for the perturbed polyharmonic operator $\mathcal{L}_{A,q}$ is to determine $A$ and $q$ in $\Omega$ from the knowledge of the Dirichlet to Neumann map $\mathcal{N}_{A,q}$.
	
	Before we proceed let us fix some notations. Here and in what follows, $\mathcal{E}'(\bar \Omega) = \{u\in \mathcal{D}'(\mathbb{R}^n: supp(u) \subseteq \bar \Omega\}$ and $W^{s,p}(\mathbb{R}^n)$ is the standard $L^p$ based Sobolev space on $\mathbb{R}^n$, $s \in \mathbb{R}$ and $1<p<\infty$, which is defined via the Bessel potential operator.  We can also define the analogous spaces $W^{s,p}(\Omega)$ for $\Omega$ a bounded open set with smooth boundary. We refer the reader to \cite{Adams} for properties of these spaces.
	
	The study of inverse problems for such first order perturbations of the polyharmonic operator was initiated in \cite{KLU}. The authors tackled the question of unique recovery of $A$  and $q$ from the knowledge of the Dirichlet-to-Neumann map. More precisely, they show that for $m \geq 2$, the set of Cauchy data \\$C_{A,q} = \{(\gamma u, \tilde \gamma u ): u\in H^{2m}(\Omega)\;\text{ with } \mathcal{L}_{A,q}u = 0\}$ determines $A$ and $q$ uniquely provided $A \in W^{1,\infty}(\Omega, \mathbb{C}^n) \cap \mathcal{E}'(\bar \Omega, \mathbb{C}^n)$ and $q \in L^{\infty}(\Omega)$. Regularities of $A$ and $q$ were substantially relaxed by the first author in \cite{Yernat} to $A\in W^{-\frac{m-2}{2},\frac{2n}{m}}(\mathbb R^n)\cap \mathcal E'(\bar \Omega)$ and $q\in W^{-\frac{m}{2}+\delta,\frac{2n}{m}}(\mathbb R^n)\cap \mathcal E'(\bar \Omega)$, $0<\delta<1/2$, for the case $m<n$. A natural question that remained open was the problem  of uniqueness in this inverse problem when the regularity of the coefficients is significantly lowered. 
	
	In this paper we successfully tackle the above question and  improve the results of \cite{KLU} and \cite{Yernat} in several directions. We show that the restriction $m<n$ in \cite{Yernat} is not necessary and that the uniqueness can in fact be proved for any $n \geq 3$ and any $m \geq 2$. Second, we substantially relax the regularity and integrability conditions for $A$ and $q$ and prove the uniqueness result for the stated inverse problem for a much broader class of coefficients. Third, we show how careful book-keeping in fact improves the result in \cite{Yernat} for $m<n$. Along the way we also reason how the class of coefficients for which uniqueness in this inverse problem can be answered using this technique, is as broad as possible and cannot be further improved. 
	
	Let us remark that the problem considered in this paper can be considered as a generalization of Calder\'{o}n's inverse conductivity problem \cite{Calderon}, also known as electrical impedance tomography, for which the question of reducing regularity has been studied extensively. In the fundamental paper by Sylvester and Uhlmann \cite{SU} it was shown that $C^2$ conductivities can be uniquely determined from boundary measurements. Successive papers have focused on weakening the regularity for the conductivity; see \cite{Brown, Paivarinta, Greenleaf, Boaz, Haberman, Pedro} for more details.
	
	As was observed in \cite{ZiqiSun}, for the case $m=1$ in \eqref{1.1}, there is a gauge invariance that prohibits uniqueness and therefore we can hope to recover $A$ and $q$ only modulo such a gauge transformation. It was shown in \cite{ZiqiSun} that such uniqueness modulo a gauge invariance is possible provided that $A \in W^{2,\infty}$, $q \in L^{\infty}$ and $dA$ satisfy a smallness condition. There have been many successive papers which have weakened the regularity assumptions on $A$ and $q$ for the case $m=1$.  The reader is referred to \cite{Mikko, Katya, NakamuraSunUhlmann, Tolmasky, Habermanmagnetic} for details.
	
	Inverse problems for higher order operators have been considered in \cite{KLU, KrupchykUhlmann, Ghosh, Krishnan, Yernat, Yernat-Yang} where unique recovery actually becomes possible. Higher order polyharmonic operators arise in the areas of physics and geometry such as the study of the Kirchoff plate equation in the theory of elasticity, and the study of the Paneitz-Branson operator in conformal geometry; for more details see ~\cite[Chapter 1]{Gazzola}.

	\subsection{Statement of Result}
	
	Throughout this paper we assume $m \geq  2$ and $n \geq 3$.
	Suppose that the first order perturbation $A$ be in $W^{-\frac{m}{2}+1,p'}(\mathbb{R}^n) \cap\, \mathcal{E}'(\bar \Omega)$, where $p'$ satisfies
	\begin{equation}\label{defition of A}
	\begin{cases}
	p' \in [2n/m, \infty) &\text{if}\quad m<n,\\
	p' \in (2,\infty) &\text{if}\quad m=n\quad \mbox{or}\quad m=n+2,\\
	p' \in [2,\infty) &\text{otherwise}.
	\end{cases}
	\end{equation}
	For a fixed $\delta$ with $0 < \delta < \frac{1}{2}$, suppose that the zeroth order perturbation $q$ be in $W^{-\frac{m}{2} + \delta,r'}(\mathbb{R}^n) \cap\, \mathcal{E}'(\bar \Omega)$, where $r'$ satisfies
	\begin{equation}\label{defition of q}
	\begin{cases}
	r' \in [2n/(m - 2 \delta), \infty), &\text{if}\quad m<n,\\
	r' \in [2n/(m - 2 \delta), \infty), &\text{if}\quad m=n,\\
	r' \in [2, \infty), &\text{if}\quad m \geq n+1.
	\end{cases}
	\end{equation}

	Before stating the main result, we consider the  bi-linear forms $B_{A}$ and $b_{q}$ on $H^m(\Omega)$ which are defined by 
	\begin{equation}\label{bilinearformdefn}
	B_A(u,v): = B_{A}^{\mathbb{R}^n}(\tilde u, \tilde v):=\langle A, \tilde v D \tilde u \rangle,\quad b_q(u,v): = b_{q}^{\mathbb{R}^n}(\tilde u, \tilde v):=\langle q, \tilde u \tilde v \rangle
	\end{equation}
	for all $u, v \in H^{m}(\Omega)$, where $\langle\cdot,\cdot\rangle$ denotes the distributional duality on $\mathbb{R}^n$ such that $\langle\cdot,\bar{\text{$\cdot$}}\rangle$ naturally extends $L^2(\mathbb R^n)$-inner product, and $\tilde u, \tilde v \in H^{m}(\mathbb{R}^n)$ are any extensions of $u$ and $v$,  respectively. In Appendix A, we show that these definitions are well defined i.e. independent of the choice of extensions $\tilde u, \tilde v$.
	Using a property of multiplication of functions in Sobolev spaces, we show that the forms $B_{A}$ and $b_{q}$ are bounded on $H^{m}(\Omega)$. We also adopt the convention that for any $z>1$, the number $z'$ is defined by $z'=z/(z-1)$.
	
	Consider the operator $D_{A}$, which is formally $A\cdot D$ where $D_{j} = -i \partial_{j}$, and the operator $m_{q}$ of multiplication by $q$. To be precise, for $u \in H^{m}(\Omega)$, $D_A(u)$ and $m_{q}(u)$ are defined as 
	$$
	\langle D_A(u),\psi \rangle_{\Omega} = B_A(u, \psi) \quad and \quad \langle m_q(u),\psi \rangle_{\Omega} = b_q(u,\psi), \quad \psi \in C^{\infty}_{0}(\Omega),
	$$
	where $\langle \cdot,\cdot \rangle_{\Omega}$ is the distribution duality on $\Omega$ such that $\langle \cdot,\bar\cdot \rangle_{\Omega}$ naturally extends $L^2(\Omega)$-inner product. The operators $D_A$ and $m_q$ are shown in Appendix A to be bounded $H^m(\Omega) \to H^{-m}(\Omega)$  and hence, standard arguments show that the operator $\mathcal{L}_{A,q} = (- \Delta)^m + D_A + m_q: H^m(\Omega) \to H^{-m}(\Omega) = (H^m_{0}(\Omega))'$ is a Fredholm operator with index $0$.
	
	For $f=(f_0,f_1,...,f_{m-1}) \in \prod_{j=0}^{m-1} H^{m - j - 1/2}(\partial \Omega)$, consider the Dirichlet problem (\ref{1.1}). If $0$ is not in the spectrum of $\mathcal{L}_{A,q}$, it is shown in Appendix B (in proof of Proposition \ref{well-posedness}) that the Dirichlet problem (\ref{1.1})  has a unique solution $u \in H^{m}(\Omega)$. We define the Dirichlet-to-Neumann map $\mathcal{N}_{A,q}$ weakly as follows
	\begin{equation}
	\langle N_{A,q}f, \bar h \rangle_{\partial \Omega} = \sum_{|\alpha| = m} \frac{m !}{\alpha !} (D^\alpha u, D^{\alpha} v_h)_{L^{2}(\Omega)} + B_A(u,\bar v_h) + b_{q}(u, \bar v_h),
	\end{equation}
	where $h = (h_0,h_1,...,h_{m-1}) \in H^{m - j - 1/2}(\partial \Omega)$, $v_h \in H^{m}(\Omega)$ is any extension of $h$ so that $\gamma v_{h} = h$, and where $\langle \cdot,\cdot \rangle_{\partial \Omega}$ is the distribution duality on $\partial \Omega$ such that $\langle \cdot,\bar\cdot \rangle_{\partial \Omega}$ naturally extends $L^2(\partial \Omega)$-inner product. It is shown in Proposition \ref{well-posedness} in Appendix~B that $\mathcal{N}_{A,q}$ is a well-defined bounded operator mapping
	$$
	\prod_{j=0}^{m-1} H^{m-j-1/2}(\partial \Omega) \to \left ( \prod_{j=0}^{m-1} H^{m-j-1/2}(\partial \Omega) \right )' = \prod_{j=0}^{m-1} H^{-m+j+1/2}(\partial \Omega).
	$$

	Our main result is as follows.
	\begin{theorem}\label{Main Result}
		Let $\Omega \subset \mathbb{R}^{n}$, $n \geq 3$ be a bounded open set with $C^{\infty}$ boundary, and let $m \geq 2$ be an integer. Let $0 <\delta<1/2$. Suppose that $A_1$, $A_2$ satisfy \eqref{defition of A} and $q_1$, $q_2$ satisfy \eqref{defition of q} and $0$ is not in the spectrums of $\mathcal{L}_{A_1,q_1}$ and $\mathcal{L}_{A_2,q_2}$. If $\mathcal{N}_{A_1,q_1} = \mathcal{N}_{A_2,q_2}$, then $A_1= A_2$ and $q_1 = q_2$. 
	\end{theorem}
	
	Detailed explanation for the assumption $\delta>0$ is given in Remark~\ref{remark::why we need delta positive}.
	
	The proof of Theorem \ref{Main Result} is structured similarly as in \cite{Yernat}. The key ingredient in the proof of Theorem \ref{Main Result} is the construction of complex geometric optics solutions for the operator $\mathcal{L}_{A,q}$ with correct decay for the remainder term. We use the method of Carleman estimates which is based on the corresponding Carleman estimates for the Laplacian, with a gain of two derivatives, due to Salo and Tzou \cite{Salo} and chain it with Proposition \ref{MultiplicationTheorem1}, which gives property of products of functions in various Sobolev spaces; to eventually obtain the desired decay.

	The idea of constructing such complex geometric optics solutions to elliptic operators goes back to the fundamental paper by Sylvester and Uhlmann \cite{SU} and has been extensively used to show unique recovery of coefficients in many inverse problems. 

	The rest of the paper is organized as follows. In Section 2 we  construct complex geometrical optics solutions for the perturbed polyharmonic operator $\mathcal{L}_{A,q}$ with $A$ and $q$ as defined in (\ref{defition of A}) and (\ref{defition of q}) respectively. This is done by deriving Carleman estimates for $\mathcal{L}_{A,q}$. Section 3 is devoted to deriving an integral identity. The proof of Theorem \ref{Main Result} is given in Section 4.  In Appendix A we study mapping properties of $D_{A}$ and $m_{q}$.  Appendix B is devoted to the well-posedness of the Dirichlet problem $\mathcal{L}_{A,q}$ with $A$ satisfying \eqref{defition of A} and $q$ satisfying \eqref{defition of q}. In Appendix C we specify why we use  Bessel potential to define fractional Sobolev spaces.
	
	
	
	
	\section{Carleman estimate and CGO solutions}\label{previous work}
	
	As a first step, we will derive Carleman estimates for the operator $\mathcal{L}_{A,q}$. We first recall the Carleman estimates for the semi-classical Laplace operator $-h^2 \Delta$ with a gain of two derivatives, as established in \cite[Lemma 2.1]{Salo}. Let $\tilde \Omega$ be an open set in $\mathbb{R}^n$ such that $\bar \Omega \subset \subset \tilde \Omega$ and let $\phi \in C^{\infty}(\tilde \Omega, \mathbb{R})$. Consider the conjugated operator $P_{\phi} = e^{\phi/h}(-h^2 \Delta)e^{-\phi/h}$ and its semi classical principal symbol
	$p_{\phi}(x, \xi) = \xi^2 + 2i \nabla \phi \cdot \xi - |\nabla \phi|^2$, $x \in \tilde \Omega$, $\xi \in \mathbb{R}^n$. Following \cite{KenigUhlmann} we say that $\phi$ is a limiting Carleman weight for $-h^2 \Delta$ in $\tilde \Omega$, if $\nabla \phi \neq 0$ in $\tilde \Omega$ and the Poisson bracket of $Re\, p_{\phi}$ and $Im\, p_{\phi}$ satisfies $\{Re\, p_\phi, Im \,p_\phi \}(x, \xi) = 0$ when $p_{\phi}(x,\xi)=0$, $(x,\xi) \in \tilde \Omega\times \mathbb{R}^n$.
	
	Before we state the Carleman estimates in \cite[Lemma 2.1]{Salo}, we define the semi-classical Sobolev norms on $\mathbb R^n$
	$$
	||u||_{H^s_{scl}(\mathbb{R}^n)}:=||\langle hD \rangle ^s u||_{L^2(\mathbb{R}^n)},
	$$
	
	where $\langle \xi \rangle = (1 + |\xi|^2)^{1/2}$ and $s \in \mathbb{R}$. 
	
	\begin{proposition}\label{mikkoprop}
		Let $\phi$ be a limiting Carleman weight for $-h^2 \Delta$ in $\tilde \Omega$ and let $\phi_{\epsilon} = \phi + \frac{h}{2 \epsilon} \phi^2$. Then for $ 0 < h \ll \epsilon \ll 1$ and $s \in \mathbb{R}$, we have 
		$$
		\frac{h}{\sqrt{\epsilon}} ||u||_{H^{s + 2}_{scl}(\mathbb{R}^n)} \leq C ||e^{\phi_{\epsilon}/h}(-h^2 \Delta) e^{-\phi_{\epsilon}/h} u||_{H^{s}_{scl}(\mathbb{R}^n)},\quad  C >0
		$$
		for all $u \in C_{0}^{\infty}(\Omega)$.
	\end{proposition}
	
	We now state a theorem on products of functions in Sobolev spaces (see Theorem~1 and Theorem~2 in \cite[Section 4.4.4]{RunstSickel}), which are used to prove  Carleman estimates stated in Proposition \ref{MultiplicationTheorem1}.
	
	\begin{proposition}\label{MultiplicationTheorem1}
		Let $0<s_1 \leq s_2$. Suppose
		\begin{enumerate}
			\item[(a)] $p^{-1} \leq p_1^{-1} + p_2^{-1} \leq 1$;
			\item[(b)] either
			$$
			\frac{n}{p} - s_1 > \begin{cases} 
			(\frac{n}{p_1} - s_1)^{+} +   (\frac{n}{p_1} - s_2)^{+} & \mbox{ if } max_{i}(\frac{n}{p_i} - s_i) >0,\\
			max_{i}(\frac{n}{p_i} -s_i) & \mbox{ otherwise } \\
			\end{cases}
			$$
			or
			\begin{align*}
			&\frac{n}{p} - s_1 = \begin{cases} 
			(\frac{n}{p_1} - s_1)^{+} +   (\frac{n}{p_1} - s_2)^{+} & \mbox{ if } max_{i}(\frac{n}{p_i} - s_1) >0,\\
			max_{i}(\frac{n}{p_i} -s_1) & \mbox{ otherwise } \\
			\end{cases}\\
			&\{i \in \{1,2\}: s_i = n/p_i \mbox{ and }p_i >1 \} = \varnothing.
			\end{align*}
			
		\end{enumerate}
		If $u \in W^{s_1,p_1}(\mathbb{R}^n)$ and $v \in W^{s_2,p_2}(\mathbb{R}^n)$, then $uv \in W^{s_1,p}(\mathbb{R}^n)$. Moreover, the point-wise multiplication of functions is a continuous bi-linear map\\
		$W^{s_1,p_1}(\mathbb{R}^n) \cdot W^{s_2,p_2}(\mathbb{R}^n) \hookrightarrow W^{s_1,p}(\mathbb{R}^n)$
		with
		\begin{equation}\label{xkcd}
		||uv||_{W^{s_1,p}(\mathbb{R}^n)} \leq C ||u||_{W^{s_1,p_1}(\mathbb{R}^n)} ||v||_{W^{s_2,p_2}(\mathbb{R}^n)}
		\end{equation}
		where the constant $C$ depends only on the various indices. 
	\end{proposition}

	With all the preliminaries behind us, we now derive Carleman estimate for the perturbed operator $\mathcal{L}_{A,q}$  when $A$ and $q$ are as in \eqref{defition of A} and \eqref{defition of q} respectively.  We have the following estimate.
	\begin{proposition}
		Let $\phi$ be a limiting Carleman weight for $-h^2 \Delta$ in $\tilde \Omega$ and suppose 
		$A$ and $q$ satisfy (\ref{defition of A}) and (\ref{defition of q}), respectively. Then for $0<h\ll 1$ ,we have
		\begin{equation}\label{Carleman estimate}
		||u||_{H^{m/2}_{scl}(\mathbb{R}^n)} \lesssim \frac{1}{h^m} ||e^{\phi/h}(h^{2m} \mathcal{L}_{A,q}) e^{-\phi/h}u||_{H^{-3m/2}_{scl}(\mathbb{R}^n)},
		\end{equation}
		for all $u \in C^{\infty}_{0}(\Omega)$.
	\end{proposition}

	\begin{proof}
		Iterate the Carleman estimate in Proposition \ref{mikkoprop} $m$ times with $s=-3m/2$ and a fixed $\epsilon >0$ sufficiently small and independent of $h$ to get the estimate
		\begin{equation} \label{ polyharmonic estimate}
		\frac{h^m}{\epsilon^{m/2}}||u||_{H^{m/2}_{scl}(\mathbb{R}^n)} \leq C || e^{\phi_{\epsilon}/h}(-h^2 \Delta)^m  e^{-\phi_{\epsilon}/h}u||_{H^{-3m/2}_{scl}(\mathbb{R}^n)},
		\end{equation}
		for all $u \in C^{\infty}_{0}(\Omega)$ and $0<h \ll \epsilon\ll 1$. Let $\psi \in C^{\infty}_{0}(\mathbb{R}^n)$ be fixed.
		
	 Let us first estimate the term involving the zeroth order perturbation $q$. By duality and Proposition \ref{MultiplicationTheorem1}, we have for any $m \geq 2$,
		\begin{equation}\label{adjointqestimate}
		\begin{aligned}
		&| \langle e^{\phi_{\epsilon}/h} h^{2m} m_{q}(e^{-\phi_{\epsilon}/h}u), \psi \rangle | \\&\leq h^{2m} ||q||_{W^{-\frac{m}{2}+\delta,r'}(\mathbb{R}^n)}||u\psi||_{W^{\frac{m}{2}-\delta,r}(\mathbb{R}^n)}\\
		&\leq C h^{2m} ||q||_{W^{-\frac{m}{2}+\delta,r'}(\mathbb{R}^n)} ||u||_{H^{\frac{m}{2} - \delta}(\mathbb{R}^n)} ||\psi||_{H^{\frac{m}{2} - \delta}(\mathbb{R}^n)} \\
		&\leq C h^{m + 2\delta}  ||u||_{H^{m/2 - \delta}_{scl}(\mathbb{R}^n)}  ||\psi||_{H^{m/2 - \delta}_{scl}(\mathbb{R}^n)} \mbox{(By definition of $H^s_{scl}(\mathbb{R}^n)$)}\\
		&\leq C h^{m}  ||u||_{H^{m/2}_{scl}(\mathbb{R}^n)}  ||\psi||_{H^{3m/2}_{scl}(\mathbb{R}^n)} \;\mbox{(Since $0<h<1$ and $\delta > 0$).}
		\end{aligned}
		\end{equation}
		
		\begin{remark}\label{rmk::remark24}
			{\it The second inequality in \eqref{adjointqestimate} follows from Proposition~\ref{MultiplicationTheorem1}. Let us break down how. 	
			\begin{itemize}
			\item For $m<n$, in Proposition~\ref{MultiplicationTheorem1}, we choose $p_1 = p_2 = 2$, $r=p \in (1,\frac{2n}{2n-m + \delta}]$, $s_1 = s_2 = \frac{m}{2} - \delta$. 
			\item For $m=n$, in Proposition~\ref{MultiplicationTheorem1} we choose $p_1 = p_2 = 2$, $r=p \in (1,\frac{2n}{n+\delta}]$, $s_1 = s_2 = \frac{m}{2} - \delta$. 
			\item Finally, for $m>n$, in Proposition~\ref{MultiplicationTheorem1} we choose $p_1 = p_2 = 2$, $r=p \in (1,2]$, $s_1 = s_2 = \frac{m}{2} - \delta$.  
			\end{itemize}
			In all the 3 cases ($m<n$, $m=n$ and $m>n$), since $q \in W^{-\frac{m}{2} + \delta,r'}(\mathbb{R}^n) \cap\, \mathcal{E}'(\bar \Omega)$ with $r'$ satisfying \eqref{defition of q}, the above choices of $p_1, p_2, p, s_1, s_2$ are justified.}
			\end{remark}
		Thus by definition of dual norm,
		\begin{equation} \label{q estimate}
		|| e^{\phi_{\epsilon}/h} h^{2m} m_{q}(e^{-\phi_{\epsilon}/h}u)||_{H^{-3m/2}_{scl}(\mathbb{R}^n)} \leq C  h^{m} ||q||_{W^{-\frac{m}{2}+\delta,r'}(\mathbb{R}^n)} ||u||_{H^{m/2}_{scl}(\mathbb{R}^n)}.
		\end{equation}

		Let us now turn our attention to the terms involving the first order perturbation $A$. For $m >2$, by duality, we have 
		\begin{align*}
		&| \langle e^{\phi_{\epsilon}/h} h^{2m} D_{A}(e^{-\phi_{\epsilon}/h}u), \psi \rangle |\\ &= | \langle h^{2m} A ,e^{\phi_{\epsilon}/h} \psi D(e^{-\phi_{\epsilon}/h}u) \rangle |  \\
		&\leq | \langle h^{2m-1}A, \psi [-u(1 + h \phi/ \epsilon) D \phi + h Du ] \rangle  \\
		&\leq C h^{2m-1} ||A||_{W^{-\frac{m}{2} +1 , p'}(\mathbb{R}^n)}\,||-u(1 + h \phi /  \epsilon) D \phi \psi + h Du \psi || _{W^{\frac{m}{2}-1,p}(\mathbb{R}^n)}.
		\end{align*}
		
		Using Proposition \ref{MultiplicationTheorem1}, we have
		\begin{align*}
		&||-u(1 + h \phi /  \epsilon) D \phi \psi + h Du \psi || _{W^{\frac{m}{2}-1,p}(\mathbb{R}^n)}\\&\leq C || -u(1 + h \phi /  \epsilon) D \phi + h Du||_{H^{\frac{m-2}{2}}(\mathbb{R}^n)}|| \psi||_{H^{\frac{m}{2}}(\mathbb{R}^n)} \\ 
		&\leq C || -u(1 + h \phi /  \epsilon) D \phi + h Du||_{H^{\frac{m-2}{2}}(\mathbb{R}^n)}|| \psi||_{H^{m/2}(\mathbb{R}^n)} \\
		&\leq Ch^{-m+1}||u||_{H^{\frac{m}{2}}_{scl}(\mathbb{R}^n)} ||\psi||_{H^{\frac{m}{2}}_{scl}(\mathbb{R}^n)} \\
		&\leq Ch^{-m+1} ||u||_{H^{\frac{m}{2}}_{scl}(\mathbb{R}^n)} ||\psi||_{H^{\frac{3m}{2}}_{scl}(\mathbb{R}^n)}.
		\end{align*}
		
		For $m=2$, we get
		\begin{align*}
	&| \langle e^{\phi_{\epsilon}/h} h^{2m} D_{A}(e^{-\phi_{\epsilon}/h}u), \psi \rangle |\\&= | \langle h^{2m} A ,e^{\phi_{\epsilon}/h} \psi D(e^{-\phi_{\epsilon}/h}u) \rangle |  \\
		&\leq | \langle h^{2m-1}A, \psi [-u(1 + h \phi/ \epsilon) D \phi + h Du ] \rangle| \\
		&\leq C h^{2m-1} ||A||_{L^{n}(\mathbb{R}^n)}||-u(1 + h \phi /  \epsilon) D \phi \psi + h Du \psi || _{L^{n'}(\Omega)}.
		\end{align*}
		
		Now, we use H\"older's inequality and Sobolev Embedding Theorem to get 
		\begin{align*}
		&| \langle e^{\phi_{\epsilon}/h} h^{2m} D_{A}(e^{-\phi_{\epsilon}/h}u), \psi \rangle |\\&\leq C  h^{2m-1} ||A||_{L^{n}(\mathbb{R}^n)}||-u(1 + h \phi /  \epsilon) D \phi + h Du  || _{L^{2}(\mathbb{R}^n)} ||\psi||_{L^{\frac{2n}{n-2}}(\mathbb{R}^n)}  \\
		&\leq C  h^{2m-1} ||A||_{L^{n}(\mathbb{R}^n)}|| -u(1 + h \phi /  \epsilon) D \phi + h Du || _{H^{\frac{m-2}{2}}(\mathbb{R}^n)} ||\psi||_{H^1(\mathbb{R}^n)}  \\
		&\leq Ch^{m} ||u||_{H^{m/2}_{scl}(\mathbb{R}^n)} ||\psi||_{H^{m/2}_{scl}(\mathbb{R}^n)} \leq Ch^m ||u||_{H^{m/2}_{scl}(\mathbb{R}^n)} ||\psi||_{H^{3m/2}_{scl}(\mathbb{R}^n)}.
		\end{align*}
		
		Thus, we can conclude that, for any $m \geq 2$, by definition of dual norm,
		$$
		||e^{\phi_{\epsilon}/h} h^{2m} D_{A}(e^{-\phi_{\epsilon}/h}u)||_{H^{-3m/2}_{scl}(\mathbb{R}^n)} \leq C h^{m} ||A||_{W^{-\frac{m}{2}+1,p'}(\mathbb{R}^n)} ||u||_{H^{m/2}_{scl}(\mathbb{R}^n)}.
		$$
		
		Combining this together with (\ref{ polyharmonic estimate}) and (\ref{q estimate}), for small enough $h >0$ and $m \geq 2$, we get
		\begin{equation}\label{Aestimate}
		||u||_{H^{m/2}_{scl}(\mathbb{R}^n)} \lesssim \frac{1}{h^m} ||e^{\phi_{\epsilon}/h}(h^{2m} \mathcal{L}_{A,q}) e^{-\phi_{\epsilon}/h} u||_{H^{-3m/2}_{scl}(\mathbb{R}^n)}.
		\end{equation}
		
		Since $e^{-\phi_{\epsilon}/h} u = e^{-\phi/h} e^{-\phi^2/ 2\epsilon} u $ and $\phi$ is smooth, we obtain (\ref{Carleman estimate}).
	\end{proof}
	
	\begin{remark}
		{\it Note that the Carleman estimate in Proposition \ref{mikkoprop} is valid for any $\tilde t \in \mathbb{R}$. We have in particular chosen $s=-3m/2$ so that $s+2m = m/2$.  The main motivation for choosing this particular value of $s$ is to get bounds on $H^{m/2}_{scl}(\mathbb{R}^n)$ norm of $u$. Though the direct problem has a solution in $H^m(\Omega)$ we only need Carleman estimates in the $H^{m/2}_{scl}$ norm.
		
		A natural question would be why in particular has $s$ been chosen so that $s+2m = m/2$.  If we choose $s+2m <m/2$ then we will have to take more regular $A$ and $q$ to ensure that we have the correct decay essentially as dictated by the hypotheses in Proposition \ref{MultiplicationTheorem1}. If we choose $s +2m >m/2$ we can no longer ensure a decay of at least $\mathcal{O}(h^m)$ for $||e^{\phi/h}(h^{2m} \mathcal{L}_{A,q}) e^{-\phi/h}u||_{H^{s}_{scl}(\mathbb{R}^n)}$ which is crucially used in the construction of complex geometric optics solutions.}
	\end{remark}

	We now use the above proved Carleman estimate to first establish an existence and uniqueness result for the inhomogeneous partial differential equation. Let $\phi \in C^{\infty}(\tilde \Omega, \mathbb{R})$ be a limiting Carleman weight for $-h^2 \Delta$. Set 
	$$
	\mathcal{L}_{\phi}:=e^{\phi/h}(h^{2m} \mathcal{L}_{A,q}) e^{-\phi/h}.
	$$
	
	By Proposition (\ref{Adjoint}), we have 
	$$
	\langle \mathcal L_{\phi}u, \bar v \rangle_{\Omega} = \langle u,\overline{\mathcal L_{\phi}^{*}v}\rangle _{\Omega}, \quad u, v \in C^{\infty}_{0}(\Omega),
	$$
	where $\mathcal{L}_{\phi}^{*} = e^{-\phi/h}(h^{2m}\mathcal{L}_{\bar A, \bar q + D \cdot \bar A})e^{\phi/h}$ is the formal adjoint of $\mathcal{L_{\phi}}$. The Carleman estimate for the first order coefficient of the adjoint operator $\mathcal{L}_{\phi}^{*}$ is the same as \eqref{Aestimate} since $\bar A$ lies in the same class as $A$. Note that the zeroth order coefficient of the adjoint operator $\mathcal{L}^{*}_{A,q}$ comprises of two terms $\bar q $ and $D \cdot \bar A$. The Carleman estimate for $\bar q $ is the same as \eqref{adjointqestimate} as $\bar q $ lies in the same class as $q$. 
	
	However $D \cdot \bar A \in W^{-\frac{m}{2}, p'}(\mathbb{R}^n) \cap \mathcal{E}'(\bar \Omega)$ where 
	$p' \geq \frac{2n}{m}$ if $m<n$, $p'>2$ if $m=n$ or $m=n+2$ and $p' \geq 2$ otherwise. The analogue of \eqref{adjointqestimate} for $D \cdot \bar A$ is as follows. 
	
	We have for $m \geq 2$, 
	\begin{equation}\label{adjointqestimateconjugate}
	\begin{aligned}
	&| \langle e^{\phi_{\epsilon}/h} h^{2m} m_{D \cdot \bar A}(e^{-\phi_{\epsilon}/h}u), \psi \rangle | \\&\leq h^{2m} ||D \cdot \bar A||_{W^{-\frac{m}{2},p'}(\mathbb{R}^n)}||u\psi||_{W^{\frac{m}{2},p}(\mathbb{R}^n)}\\
	&\leq C h^{2m} ||D \cdot \bar A||_{W^{-\frac{m}{2},p'}(\mathbb{R}^n)} ||u||_{H^{\frac{m}{2}}(\mathbb{R}^n)} ||\psi||_{H^{\frac{m}{2}}(\mathbb{R}^n)} \\
	&\leq C h^{m}  ||u||_{H^{m/2}_{scl}(\mathbb{R}^n)}  ||\psi||_{H^{m/2}_{scl}(\mathbb{R}^n)} \mbox{(By definition of $H^s_{scl}(\mathbb{R}^n)$)}\\
	&\leq C h^{m}  ||u||_{H^{m/2}_{scl}(\mathbb{R}^n)}  ||\psi||_{H^{3m/2}_{scl}(\mathbb{R}^n)}
	\end{aligned}
	\end{equation}
		\begin{remark}\label{rmk::remark24}
		{\it The second inequality in \eqref{adjointqestimateconjugate} follows from Proposition~\ref{MultiplicationTheorem1}. Let us break down how. 	
			\begin{itemize}
				\item For $m<n$, in Proposition~\ref{MultiplicationTheorem1}, we choose $p_1 = p_2 = 2$, $p \in (1,\frac{2n}{2n-m}]$, $s_1 = s_2 = \frac{m}{2}$. 
				\item For $m=n$ and $m=n+2$, in Proposition~\ref{MultiplicationTheorem1} we choose $p_1 = p_2 = 2$, $p \in (1,2)$, $s_1 = s_2 = \frac{m}{2}$. 
				\item Finally, for $m = n+1$ or otherwise, in Proposition~\ref{MultiplicationTheorem1} we choose $p_1 = p_2 = 2$, $p \in (1,2]$, $s_1 = s_2 = \frac{m}{2}$.  
			\end{itemize}
			In all the 3 cases ($m<n$, $m=n$ or $m=n_2$ and $m>n$), since $D \cdot \bar A \in W^{-\frac{m}{2},p'}(\mathbb{R}^n) \cap\, \mathcal{E}'(\bar \Omega)$ with $p'$ satisfying \eqref{defition of A}, the above choices of $p_1, p_2, p, s_1, s_2$ are justified.}
	\end{remark}
	Let us now convert the Carleman estimate (\ref{Carleman estimate}) for $\mathcal{L}_{\phi}^{*}$ into a solvability result for $\mathcal{L}_{\phi}$. For $s \geq 0$, we define semi-classical Sobolev norms on a smooth bounded domain $\Omega$ as
	\begin{align*}
	||u||_{H^s_{scl}(\Omega)} &:= \inf_{v \in H^s_{scl}(\mathbb{R}^n), \; v|_{\Omega}=u} \;||v||_{H^s_{scl}(\mathbb{R}^n)},\\
	||u||_{H^{-s}_{scl}(\Omega)} &:= \sup_{ 0\neq \phi \in C^{\infty}_{0}(\Omega)} \frac{|\langle u ,\bar \phi \rangle|}{||\phi||_{H^{s}_{scl}(\Omega)}}.
	\end{align*}
	
	
	\begin{proposition}\label{solvability result}
		Let $A$ and $q$ satisfy the conditions in (\ref{defition of A}) and (\ref{defition of q}) respectively and let $\phi$ be a limiting Carleman weight for $-h^2 \Delta$ on $\tilde \Omega$. If $ h>0$ is small enough, then for any $v \in H^{-\frac{m}{2}}(\Omega)$, there is a solution $u \in H^{\frac{m}{2}}(\Omega)$ of the equation
		$$
		e^{\phi/h}(h^{2m} \mathcal{L}_{A,q})e^{-\phi/h} u = v \quad\text{in}\quad \Omega,
		$$
		which satisfies 
		$$
		||u||_{H^{\frac{m}{2}}_{scl}(\Omega)} \lesssim \frac{1}{h^m} ||v||_{H^{-\frac{m}{2}}_{scl}(\Omega)}.
		$$
	\end{proposition}
	\begin{proof}
		Let $D = \mathcal{L}_{\phi}^{*}(C_{0}^{\infty}(\Omega))$ and consider the linear functional\\ $L: D \to \mathbb{C}$, $L(\mathcal{L}_{\phi}^{*}w) = \langle w,  v\rangle_{\Omega}$ for $w \in C^{\infty}_{0}(\Omega)$. 
		By the Carleman estimate (\ref{Carleman estimate}), 
		\begin{align*}
		|L(\mathcal{L}_{\phi}^{*} w) |&\leq ||w||_{H^{m/2}_{scl}(\mathbb{R}^n)} ||v||_{H^{-m/2}_{scl}(\Omega)}\\& \leq \frac{C}{h^m}  ||\mathcal{L}_{\phi}^{*} w||_{H^{-m/2}_{scl}(\mathbb{R}^n)} ||v||_{H^{-m/2}_{scl}(\Omega)}.
		\end{align*}
		
		Hahn-Banach theorem ensures that there is a bounded linear functional \\$\tilde L : H^{-m/2}(\mathbb{R}^n) \to \mathbb{C}$ satisfying $\tilde L = L$ on $D$ and $||\tilde L|| \leq Ch^{-m} ||v||_{H^{-m/2}_{scl}(\Omega)}$. 
		
		By the Riesz Representation theorem there is $u \in H^{m/2}(\mathbb{R}^n)$ such that for all $\psi \in H^{-m/2}(\mathbb{R}^n)$, $\tilde L (\psi) = \langle \psi,\bar u \rangle_{\mathbb{R}^n}$ and 
		$$
		||u||_{H^{m/2}_{scl}(\mathbb{R}^n)}\leq  \frac{C}{h^m} ||v||_{H^{-m/2}_{scl}(\Omega)}.
		$$
		
		Let us show $\mathcal L_{\phi}u = v $ in $\Omega$. For arbitrary $w \in C^{\infty}_{0}(\Omega)$,
		$$
		\langle \mathcal L_\phi u,\bar w \rangle_{\Omega} = \langle u ,\overline{\mathcal L^{*}_{\phi}w}\rangle_{\mathbb{R}^n} = \overline{\tilde L(\mathcal L^{*}_{\phi}w)} = \overline{L(\mathcal L^{*}_{\phi}w)} = \overline{\langle w, v\rangle}_{\Omega} = \langle v,\bar w\rangle_{\Omega}.
		$$
		
		This finishes the proof.
	\end{proof}
	We now wish to construct complex geometric optics solutions for the equation $\mathcal{L}_{A,q}u =0$ in $\Omega$ with $A$ and $q$ as defined in (\ref{defition of A}) and (\ref{defition of q}) respectively using the solvability result Proposition \ref{solvability result}. These are solutions of the form
	\begin{equation}\label{CGO form}
	u(x, \zeta;h)=e^{\frac{i x \cdot \zeta}{h}}(a(x,\zeta;h) + h^{m/2}r(x,\zeta;h)),
	\end{equation}
	where $\zeta \in \mathbb{C}^n$ is such that $\zeta \cdot \zeta = 0$, $|\zeta| \sim 1$, $a \in C^{\infty}(\bar \Omega)$ is an amplitude, $r$ is a correction term, and $h >0$ is a small parameter.
	
	Conjugating $h^{2m} \mathcal{L}_{A,q}$ by $e^{\frac{i x \cdot \zeta}{h}}$, we get
	\begin{equation} \label{conjugation by exp}
	e^{-\frac{i x \cdot \zeta}{h}} h^{2m} \mathcal{L}_{A,q} e^{\frac{i x \cdot \zeta}{h}} = (-h^2 \Delta  - 2 i \zeta \cdot h \nabla)^m + h^{2m}D_{A} + h^{2m-1}m_{A \cdot \zeta} + h^{2m}m_{q}.
	\end{equation}
	
	Following \cite{Katya}, we shall consider $\zeta$ depending slightly on $h$, i.e $\zeta = \zeta_{0} + \zeta_{1}$ with $\zeta_0$ independent of $h$ and $\zeta_1 = \mathcal{O}(h)$ as $h \to 0$. We also assume that \\$|Re\, \zeta_0| = |Im \zeta_0| = 1$. Then we can write (\ref{conjugation by exp}) as 
	\begin{align*}
	&e^{-\frac{i x \cdot \zeta}{h}} h^{2m} \mathcal{L}_{A,q} e^{\frac{i x \cdot \zeta}{h}} \\&= (-h^2 \Delta  - 2 i \zeta_0 \cdot h \nabla - 2 i \zeta_1 \cdot h \nabla)^m + h^{2m}D_A+ h^{2m-1}m_{A \cdot (\zeta_1 + \zeta_0)} + h^{2m}m_q.
	\end{align*}
	
	Observe that (\ref{CGO form}) is a solution to $\mathcal{L}_{A,q}=0$ if and only if
	$$
	e^{-\frac{i x \cdot \zeta}{h}} h^{2m} \mathcal{L}_{A,q} (e^{\frac{i x \cdot \zeta}{h}} h^{m/2}r)= - e^{-\frac{i x \cdot \zeta}{h}} h^{2m} \mathcal{L}_{A,q} (e^{\frac{i x \cdot \zeta}{h}}a),
	$$
	
	and hence if and only if
	\begin{equation}\label{*}
	\begin{aligned}
	&e^{-\frac{i x \cdot \zeta}{h}} h^{2m} \mathcal{L}_{A,q} (e^{\frac{i x \cdot \zeta}{h}} h^{m/2}r)\\
	&= - \sum_{k=0}^{m} \frac{m!}{k! (m-k)!} ( -h^2 \Delta - 2 i \zeta_1 \cdot h \nabla)^{m-k} \,(-2i \zeta_0 \cdot h \nabla)^{k}a  \\
	&\qquad -h^{2m}D_{A}a - h^{2m-1}m_{A \cdot (\zeta_0 + \zeta_1)} a - h^{2m}m_q a.
	\end{aligned}
	\end{equation}
	
	Our goal is get a decay of at least $\mathcal{O}(h^{m+m/2})$ in $H^{-m/2}_{scl}(\Omega)$ norm on the right-hand side of \eqref{*}. The terms $h^{2m}D_A a$, $h^{2m-1}m_{A \cdot (\zeta_0 + \zeta_1)}a $ and $h^{2m}m_q a$ will eventually give us a decay of $\mathcal{O}(h^{m+m/2})$ provided $m \geq 2$.

	({\it For a smooth enough first order perturbation of the polyharmonic operator we only need an $\mathcal{O}(h^{m+1})$ decay but here we need a stronger decay of $\mathcal{O}(h^{m +m/2})$ essentially because our coefficients are less regular. See  Remark \ref{rmk::remark33} for more details.)}
	
	If $a \in C^{\infty}(\bar \Omega)$ satisfies 
	$$
	(\zeta_{0} \cdot \nabla)^j a = 0 \quad \text{in}\quad\Omega
	$$
	for some $j \geq 1$, then since $\zeta_1 = \mathcal{O}(h)$, the lowest order of $h$ on the right-hand side of \eqref{*} is $j-1 + 2(m-j+1) = 2m-j+1$ provided $j \geq 2$. We will hence obtain an overall decay of $\mathcal{O}(h^{m + m/2})$ on the right-hand side of \eqref{*} provided $j \leq 1+m/2$.
	
	 Since $m \geq 2$, we choose $j=2$ to get the following transport equation,
	\begin{equation}\label{transport}
	(\zeta_{0} \cdot \nabla)^2 a = 0 \quad \text{in}\quad \Omega.
	\end{equation}
	
	Such choice of $a$ is clearly possible. We thus obtain the following equation for $r$,
	\begin{align*}
	e^{-\frac{i x \cdot \zeta}{h}} &h^{2m} \mathcal{L}_{A,q} (e^{\frac{i x \cdot \zeta}{h}} h^{m/2}r) \\
	&= - ( -h^2 \Delta - 2 i \zeta_1 \cdot h \nabla)^m a - m( -h^2 \Delta - 2 i \zeta_1 \cdot h \nabla)^{m-1} \,(-2i \zeta_0 \cdot h \nabla)a \\
	&\qquad -h^{2m}D_{A}(a) - h^{2m-1}m_{A \cdot (\zeta_0 + \zeta_1)}(a) - h^{2m}m_q(a) :=g.
	\end{align*}
	
	We complete the proof by showing $||g||_{H^{-m/2}_{scl}(\Omega)} = \mathcal{O}(h^{m +m/2})$. We will estimate each term separately. 
	
	Suppose that $ \psi \in C^{\infty}_{0}(\Omega)$ and $\psi \neq 0$. By Cauchy-Schwarz inequality and the fact that $\zeta_1 = \mathcal{O}(h)$ and $\zeta_0 = \mathcal{O}(1)$ we get
	\begin{equation}\label{226}
	\begin{aligned}
	\big| (  - ( -h^2 \Delta &- 2 i \zeta_1 \cdot h \nabla)^m a - m( -h^2 \Delta - 2 i \zeta_1 \cdot h \nabla)^{m-1} \,(-2i \zeta_0 \cdot h \nabla)a, \psi )_{L^2(\Omega)} \big|\\
	&= \mathcal{O}(h^{m+m/2}) ||\psi||_{L^2(\Omega)} = \mathcal{O}(h^{m+m/2}) ||\psi||_{H^{m/2}_{scl}(\Omega)}.
	\end{aligned}
	\end{equation}
	
	For $m>2$ we have
	\begin{equation}\label{227}
	\begin{aligned}
	&| \langle h^{2m-1} m_{A \cdot (\zeta_0 + \zeta_1)}(a), \psi  \rangle_{\Omega}| \leq C h^{2m-1} ||A||_{W^{-\frac{m}{2}+1,p'}(\mathbb{R}^n)} ||a \psi||_{W^{\frac{m}{2}-1,p}(\mathbb{R}^n)}\\
	& \leq C h^{2m-1} ||a \psi||_{W^{\frac{m}{2}-1,p}(\mathbb{R}^n)}\,\mbox{(as $A \in W^{-\frac{m}{2}+1,p'}(\mathbb{R}^n)$)} \\
	&\leq Ch^{m + m/2 + m/2 -1}||\psi||_{H^{\frac{m}{2}-1}(\mathbb{R}^n)}\,\mbox{(by Proposition \ref{MultiplicationTheorem1}) } \\
	& = \mathcal{O}(h^{m+m/2})||\psi||_{H^{m/2}_{scl}(\mathbb{R}^n)}=\mathcal{O}(h^{m+m/2})||\psi||_{H^{m/2}_{scl}(\Omega)}.
	\end{aligned}
	\end{equation}
	
	\begin{remark}
	{\it Here, we have used Proposition \ref{MultiplicationTheorem1} for $m<n$ with $p_1 = p_2 = 2$, $p \in (1,\frac{2n}{2n-m}) $, $s_1 = \frac{m}{2}-1 , s_2 = \frac{m}{2}$. For $m=n$, we choose $p_1 = p_2 = 2$, $p \in (1,2)$, $s_1 = \frac{m}{2}-1 , s_2 = \frac{m}{2}$. And for $m>n$, we choose $p_1 = p_2 = 2$, $p \in (1,2]$, $s_1 = \frac{m}{2}-1 , s_2 = \frac{m}{2}$. Thus \eqref{227} is justified for all $m>2$.}
\end{remark}
	Similarly, for $m>2$, we also have
	\begin{equation}\label{228}
	\begin{aligned}
	&| \langle h^{2m} D_{A}(a), \psi \rangle_{\Omega}|\leq C h^{2m} ||A||_{W^{-\frac{m}{2}+1,p'}(\mathbb{R}^n)} || \psi Da ||_{W^{\frac{m}{2}-1,p}(\mathbb{R}^n)}\\
	&\leq  C h^{ m+ m/2 + m/2} ||A||_{W^{-\frac{m}{2}+1,p'}(\mathbb{R}^n)} ||\psi||_{H^{m/2}(\mathbb{R}^n)}\,\mbox{(by Proposition \ref{MultiplicationTheorem1}) } \\
	&=\mathcal{O}(h^{m+m/2}) ||\psi||_{H^{m/2}_{scl}(\Omega)}.
	\end{aligned}
	\end{equation}
	
	For $m=2$, we have 
	\begin{equation}\label{229}
	\begin{aligned}
	| \langle h^{2m-1} m_{A \cdot (\zeta_0 + \zeta_1)}(a), \psi  \rangle_{\Omega}|&\leq C h^{2m-1} ||A||_{L^{n}(\mathbb{R}^n)} ||a \psi||_{L^{n'}(\Omega)}\notag\\
	& \leq C h^{2m-1} ||a \psi||_{L^{n'}(\Omega)} \\
	&\leq Ch^{2m -1}||\psi||_{L^{2}(\Omega)}\,\mbox {(by H\"older's inequality)} \\
	& = \mathcal{O}(h^{m+m/2})||\psi||_{H^{m/2}_{scl}(\Omega)}
	\end{aligned}
	\end{equation}
	
	and
	\begin{equation}\label{230}
	\begin{aligned}
	| \langle h^{2m} D_{A}(a), \psi \rangle_{\Omega}|&\leq C h^{2m} ||A||_{L^{n}(\mathbb{R}^n)} || \psi Da ||_{L^{n'}(\Omega)} \\
	&\leq  C h^{ m+ m/2 + m/2} ||A||_{L^{n}(\mathbb{R}^n)} ||\psi||_{L^2(\Omega)}\,\mbox{(by H\"older's inequality)}  \\
	&= \mathcal{O}(h^{m +m/2}) ||\psi||_{H^{m/2}_{scl}(\Omega)}.
	\end{aligned}
	\end{equation}
	
	We also have, for any $m \geq 2$,
	\begin{equation}\label{231}
	\begin{aligned}
	| \langle h^{2m} m_{q}a, \psi \rangle_{\Omega}| &\leq C h^{2m} ||q||_{W^{-\frac{m}{2} + \delta,r'}(\mathbb{R}^n)}||a \psi||_{W^{\frac{m}{2}-\delta,r}(\mathbb{R}^n)}\,\mbox{(as $q \in W^{-\frac{m}{2}+\delta,r'}(\mathbb{R}^n)$)}\\
	&\leq C h^{m + m/2 + m/2} ||q||_{W^{-\frac{m}{2} + \delta,r'}(\mathbb{R}^n)} ||\psi||_{H^{\frac{m}{2}}(\mathbb{R}^n)}\,\mbox{(by Proposition \ref{MultiplicationTheorem1}) } \\
	& = \mathcal{O}(h^{m + m/2}) ||\psi||_{H^{m/2}_{scl}(\mathbb{R}^n)}=\mathcal{O}(h^{m + m/2}) ||\psi||_{H^{m/2}_{scl}(\Omega)}.
	\end{aligned}
	\end{equation}
	
	Combining the estimates (\ref{226} - \ref{231}) we conclude that for any $m \geq 2$
	$$
	||g||_{H^{-m/2}_{scl}(\Omega)}  =\mathcal{O}(h^{m+m/2}).
	$$
	
	Using this and Proposition \ref{solvability result}, for $h>0$ small enough, we can conclude that there exists $r \in H^{m/2}(\Omega)$ solving
	$$
	e^{-\frac{i x \cdot \zeta}{h}} h^{2m} \mathcal{L}_{A,q} (e^{\frac{i x \cdot \zeta}{h}} h^{m/2}r) = -e^{-\frac{i x \cdot \zeta}{h}} h^{2m} \mathcal{L}_{A,q} (e^{\frac{i x \cdot \zeta}{h}} a),
	$$
	
	such that
	\begin{align*}
	||h^{m/2}r||_{H^{m/2}_{scl}(\Omega)} &\lesssim \frac{1}{h^m} ||e^{\frac{-i x \cdot \zeta}{h}} h^{2m} \mathcal{L}_{A,q}(e^{\frac{i x \cdot \zeta}{h}}a)||_{H^{-m/2}_{scl}(\Omega)}\\&\lesssim \frac{1}{h^m} ||g||_{H^{-m/2}_{scl}(\Omega)} = \mathcal{O}(h^{m/2}).
	\end{align*}
	
	Therefore $||r||_{H^{m/2}_{scl}(\Omega)} = \mathcal{O}(1)$. Hence we have the following result.
	
	\begin{proposition}\label{CGO}
		Let $\Omega \subset \mathbb{R}^{n}$, $n \geq 3 $ be a bounded open set with smooth boundary and let $m$ be an integer so that $m \geq 2$. Suppose $A$ and $q$ satisfy (\ref{defition of A}) and (\ref{defition of q}), respectively, and let $\zeta \in \mathbb{C}^n$ be such that $\zeta \cdot \zeta = 0$, $\zeta = \zeta_0 + \zeta_1$ with $\zeta_0$ independent of $h$ and $\zeta_1 =\mathcal{O}(h)$ as $h \to 0$. Then for all $h>0$ small enough, there exists a solution $u(x, \zeta;h) \in H^{m/2}(\Omega)$ to the equation $\mathcal{L}_{A,q}u=0$ of the form
		$$
		u(x, \zeta;h) = e^{\frac{i x \cdot \zeta}{h}}( a(x,\zeta_0) + h^{m/2} r(x,\zeta;h)),
		$$
		where $a(\cdot,\zeta_0) \in C^{\infty}(\bar \Omega)$ satisfies (\ref{transport}) and the correction term $r$ is such that\\ $||r||_{H^{m/2}_{scl}(\Omega)} = \mathcal{O}(1)$ as $h \to 0$.
	\end{proposition}
\section{Integral Identity}
We first do a standard reduction to a larger domain. For the proof we follow \cite[Proposition~3.2]{Katya}.
\begin{proposition}\label{biggerdomain}
	Let $\Omega$, $\Omega ' \subset \mathbb{R}^n$ be two bounded open sets such  that $\Omega \subset \subset \Omega '$ and $\partial \Omega$ and $\partial \Omega '$ are smooth. Let $A_1$, $A_2$ and $q_1,q_2$ satisfy (\ref{defition of A}) and (\ref{defition of q}), respectively.  If $\mathcal{N}_{A_1, q_1} = \mathcal{N}_{A_2,q_2}$, then $\mathcal{N}'_{A_1,q_1} = \mathcal{N}'_{A_2,q_2}$ where $\mathcal{N}'_{A_j,q_j}$ denotes the Dirichlet-to-Neumann map for $\mathcal{L}_{A_j,q_j}$ in $\Omega '$, $j=1,2$.  
\end{proposition}
\begin{proof}
	Let $f' \in \prod_{j=0}^{m-1} H^{m-j-1/2}(\partial \Omega)$ and let $v_{1}' \in H^{m}(\Omega ')$ be the unique solution (See Appendix B for justification of this statement) to $\mathcal{L}_{A_1,q_1}v_1' = 0$ in $\Omega'$ with $\gamma ' v_1' =f'$ on $\partial \Omega '$ where $\gamma '$ denotes the Dirichlet trace on $\partial \Omega '$. Let $v_1 = \restr{v_1'}{\Omega} \in H^{m}(\Omega)$ and let $f = \gamma v_1$. By the well-posedness result in Appendix B, we can guarantee the existence of a unique $v_2 \in H^{m}(\Omega)$ so that $\mathcal{L}_{A_2,q_2}v_2 = 0$ and $\gamma v_2 = \gamma v_1 = f$. Thus $\phi = v_2 - v_1 \in H^{m}_{0}(\Omega) \subset H^m_{0}(\Omega ')$. Define 
	$$
	v_2' = v_1' + \phi \in H^m(\Omega').
	$$
	Note that $v_2' = v_2$ in $\Omega$ and $\gamma ' v_2 ' = \gamma ' v_1' = f'$ on $\partial \Omega '$.
	
	We now show that $\mathcal{L}_{A_2, q_2} v_2' = 0$ in $\Omega '$. Let $\psi \in C^{\infty}_{0}(\Omega ')$. We then have
	$$
	\langle \mathcal{L}_{A_2,q_2} v_2',\bar \psi \rangle_{\Omega '} = \sum_{|\alpha| = m} \frac{m!}{\alpha !} (D^{\alpha} v_2', D^\alpha \psi)_{L^2(\Omega ')} + \langle D_{A_2}(v_2 '), \bar \psi \rangle_{\Omega '} + \langle m_{q_2}(v_2 '),\bar \psi \rangle_{\Omega '}.
	$$
	
	Since $A_2$ and $q_2$ are compactly supported in $\bar \Omega$ and $\phi \in H^{m}_{0}(\Omega)$, we can rewrite the above equality as
	\begin{align*}
	&\langle \mathcal{L}_{A_2,q_2} v_2',\bar \psi \rangle_{\Omega '}= \sum_{|\alpha| = m} \frac{m!}{\alpha !} (D^{\alpha} v_1', D^\alpha \psi)_{L^2(\Omega ')} + \sum_{|\alpha| = m} \frac{m!}{\alpha !} (D^{\alpha} \phi, D^\alpha (\restr{\psi}{\Omega})_{L^2(\Omega )}\\
	&\qquad\qquad + B_{A_2}(v_2',\restr{\bar \psi}{\Omega})+b_{q_{2}}(v_2',\restr{\bar \psi}{\Omega}) \\
	&=\sum_{|\alpha| = m} \frac{m!}{\alpha !} (D^{\alpha} v_1', D^\alpha \psi)_{L^2(\Omega ')} - \sum_{|\alpha| = m} \frac{m!}{\alpha !} (D^{\alpha} v_1, D^\alpha (\restr{\psi}{\Omega}))_{L^2(\Omega )} \\
	&\qquad\qquad+ \sum_{|\alpha| = m} \frac{m!}{\alpha !} (D^{\alpha} v_2, D^\alpha (\restr{\psi}{\Omega}))_{L^2(\Omega )}+ B_{A_2}(v_2',\restr{\bar \psi}{\Omega}) + b_{q_{2}}(v_2',\restr{\bar \psi}{\Omega}).
	\end{align*}
	
	Note that 
	$$
	\langle \mathcal{N}_{A_2,q_2}f, \overline{\gamma(\restr{\psi}{\Omega})} \rangle_{\partial \Omega}=\sum_{|\alpha| = m} \frac{m!}{\alpha !} (D^{\alpha} v_2, D^\alpha (\restr{\psi}{\Omega})_{L^2(\Omega )}+ B_{A_2}(v_2',\restr{\bar \psi}{\Omega}) + b_{q_{2}}(v_2',\restr{\bar \psi}{\Omega}).
	$$
	Hence, we have 
	\begin{multline*}
	\langle \mathcal{L}_{A_2,q_2} v_2',\bar \psi \rangle_{\Omega '}= \sum_{|\alpha| = m} \frac{m!}{\alpha !} (D^{\alpha} v_1', D^\alpha \psi)_{L^2(\Omega ')} -\sum_{|\alpha| = m} \frac{m!}{\alpha !} (D^{\alpha} v_1, D^\alpha (\restr{\psi}{\Omega}))_{L^2(\Omega )}\\
	+ \langle \mathcal{N}_{A_2,q_2}f, \gamma(\restr{\psi}{\Omega}) \rangle_{\partial \Omega}.
	\end{multline*}
	
	Since 
	$$
	\langle \mathcal{N}_{A_2,q_2}f, \overline{\gamma(\restr{\psi}{\Omega})} \rangle_{\partial \Omega} = \langle \mathcal{N}_{A_1,q_1}f, \gamma(\restr{\psi}{\Omega}) \rangle_{\partial \Omega}
	$$
	
	and
	\begin{align*}
	&\langle \mathcal{N}_{A_1,q_1}f, \overline{\gamma(\restr{\psi}{\Omega}) }\rangle_{\partial \Omega} \\& = \sum_{|\alpha| = m} \frac{m!}{\alpha !} (D^{\alpha} v_1, D^\alpha (\restr{\psi}{\Omega}))_{L^2(\Omega )}+ B_{A_1}(v_1, \restr{\bar \psi}{\partial \Omega})+b_{q_1}(v_1,\restr{\bar \psi}{\Omega}).
	\end{align*}
	
	We get
	$$
	\langle \mathcal{L}_{A_2,q_2} v_2',\bar \psi \rangle_{\Omega '} = \sum_{|\alpha |= m} \frac{m!}{\alpha !} (D^{\alpha} v_1', D^\alpha \psi)_{L^2(\Omega ')} + B_{A_1}(v_1, \restr{\bar \psi}{\partial \Omega})+ b_{q_1}(v_1,\restr{\bar \psi}{\Omega}).
	$$
	
	Using the fact $A_1$ and $q_1$ are compactly supported in $\bar \Omega$, we obtain
	\begin{align*}
	\langle \mathcal{L}_{A_2,q_2} v_2',\bar \psi \rangle_{\Omega '}&= \sum_{|\alpha| = m} \frac{m!}{\alpha !} (D^{\alpha} v_1', D^\alpha \psi)_{L^2(\Omega ')} + \langle D_{A_1}(v_1 '), \bar \psi \rangle_{\Omega '} + \langle m_{q_1}(v_1 '), \bar \psi \rangle_{\Omega '} \notag \\
	&=\langle \mathcal{L}_{A_1,q_1} v_1',\bar \psi \rangle_{\Omega '}=0.
	\end{align*}
	
	Using exact same arguments, one can show that $\mathcal{N}'_{A_2,q_2}f' = \mathcal{N}'_{A_1,q_1}f'$ on $\partial \Omega '$, which finishes the proof. 
\end{proof}

Proposition \ref{biggerdomain} allows us to lift the equality of Dirichlet-to-Neumann maps on to a larger domain. 
We now derive the following integral identity based on the assumption that $\mathcal{N}_{A_1,q_1} = \mathcal{N}_{A_2,q_2}$.

\begin{proposition}\label{reduction larger domain}
	Let $\Omega \subset \mathbb{R}^n$, $n \geq 3$ be a bounded open set with smooth boundary. Assume that $A_{1}, A_{2}$ and $q_1$, $q_2$ satisfy (\ref{defition of A}) and (\ref{defition of q}), respectively. If $\mathcal{N}_{A_1,q_1} = \mathcal{N}_{A_2,q_2}$, then the following integral identity holds
	$$
	\langle D_{A_2 - A_1}(u_2),\bar v \rangle_{\Omega} + \langle m_{q_2 -q_1}(u_2), \bar v\rangle_{\Omega} = 0,
	$$
	for any $u_2,v \in H^m(\Omega)$ satisfying $\mathcal{L}_{A_2,q_2}u_2 = 0$ in $\Omega$ and $\mathcal{L}^{*}_{A_1,q_1} v = 0$ in $\Omega$, respectively. Recall that $
	\mathcal{L}^{*}_{A,q} = \mathcal{L}_{\bar A, \bar q + D \cdot \bar A}$ is the formal adjoint of $\mathcal{L}_{A,q}$.
\end{proposition}

\begin{proof}
	Let $v$ satisfy $\mathcal{L}^{*}_{A_1,q_1}  v=0$ in $\Omega$. Let $u_1$ solve $\mathcal{L}_{A_1,q_1} u_1 = 0$ in $\Omega$. Since $\mathcal{N}_{A_1,q_1} = \mathcal{N}_{A_2,q_2}$, we choose $u_2 \in H^m (\Omega)$ solving $\mathcal{L}_{A_2,q_2}u_2$ so that $\gamma u_1 = \gamma u_2$ and $\mathcal{N}_{A_1,q_1} \gamma u_1 = \mathcal{N}_{A_2,q_2} \gamma u_2$. It then follows that 
	$$
	\mathcal{L}_{A_1,q_1}(u_1 - u_2) = D_{A_2 - A_1}(u_2) + m_{q_2 - q_1}(u_2).
	$$
	
	Since 
	$$
	\langle \mathcal{L}_{A,q} u , \bar v \rangle_{\Omega} = \langle u, \overline{\mathcal{L}^{*}_{A,q} v} \rangle_{\Omega},
	$$
	we get the desired identity.
\end{proof}

\section{Concluding steps}

To show $A_1 = A_2$, we will need to use Poincare lemma for currents \cite{deRham} which requires the domain to be simply connected. Therefore, we reduce the problem to larger simply connected domain, in particular to a ball.

Let us now fix $\Omega' = B$ to be an open ball in $\mathbb{R}^n$ such that $\Omega \subset \subset B$. Note that by Proposition \ref{biggerdomain}, $\mathcal{N}^{\Omega}_{A_1,q_1} = \mathcal{N}^{\Omega}_{A_2,q_2}$ implies $\mathcal{N}'^{B}_{A_1,q_1} = \mathcal{N}'^{B}_{A_2,q_2}$, where $\mathcal{N}'^{B}_{A_j,q_j}$ denotes the Dirichlet-to-Neumann map for $\mathcal{L}_{A_j,q_j}$ in $B$, $j=1,2$. 

Moreover, if  $A_{1}, A_{2}$ and $q_1$, $q_2$ satisfy (\ref{defition of A}) and (\ref{defition of q}), respectively, then $A_1$, $A_2$ and $q_1$, $q_2$ satisfy the same conditions for the larger domain $B$ too. Thus applying all the analysis up to Proposition \ref{reduction larger domain} gives us
 $\mathcal{N}'^{B}_{A_1,q_1} = \mathcal{N}'^{B}_{A_2,q_2}.$  By Proposition \ref{reduction larger domain}, the following integral identity holds
\begin{equation}\label{identity on larger domain}
B^{B}_{A_2 - A_1}(u_2, \bar v) + b^{B}_{q_2 - q_1}(u_2, \bar v) = 0,
\end{equation}
for any $u_2,v \in H^m(B)$ satisfying $\mathcal{L}_{A_2,q_2}u_2 = 0$ in $B$ and $\mathcal{L}^*_{A_1,q_1} v = 0$ in $B$, respectively. Henceforth $B_{A_2 - A_1}^{B}$ and $b^{B}_{q_2 -  q_1}$ denotes the bi-linear forms corresponding to $A_2 - A_1$ and $q_2 - q_1$ (these are shown to be well-defined for any open bounded domain $\Omega \subset \mathbb{R}^n$ with smooth boundary in Appendix A) in the ball B.

The key idea in the uniqueness result is to use complex geometric optics solutions $u_2$ to $\mathcal{L}_{A_2,q_2}u_2=0$ in $B$ and $v$ to $\mathcal{L}^{*}_{A_1,q_1} v=0$ in $B$ and plug them in the integral identity (\ref{identity on larger domain}). In order to construct these solutions, consider $\xi,\mu_1,\mu_2 \in \mathbb{R}^n$ such that $|\mu_1| = |\mu_2| = 1$ and $\mu_1 \cdot \mu_2 = \xi \cdot \mu_1 = \xi \cdot \mu_2 = 0$. For $h >0$, set
$$
\zeta_2 = \frac{h \xi}{2} + \sqrt{1 - h^2 \frac{|\xi|^2}{4}} \mu_1 + i \mu_2, \quad \zeta_1 = -\frac{h \xi}{2} + \sqrt{1 - h^2 \frac{|\xi|^2}{4}} \mu_1 - i \mu_2.
$$
Note that we have $\zeta_2 = \mu_1 + i \mu_2 + \mathcal{O}(h)$, $\zeta_1 = \mu_1 - i \mu_2 + \mathcal{O}(h)$, $\zeta_j \cdot \zeta_j = 0$, $j=1,2$ and $\zeta_2 - \bar \zeta_1 = h \xi$.

By Proposition \ref{CGO}, for all $h>0$ small enough, there are solutions $u_2(\cdot,\zeta_2;h)$ and$v(\cdot,\zeta_1;h)$ in $H^m(B)$ to the equations $\mathcal{L}_{A_2,q_2}u_2 = 0$ and $\mathcal{L}^{*}_{A_1,q_1} v = 0$ in $B$,
respectively, of the form
\begin{align*}
v(x, \zeta_1;h) &= e^{\frac{i x \cdot \zeta_1}{h}}(a_1(x,\mu_1 + i \mu_2) + h^{m/2}r_1(x,\zeta_1;h)), \\
u_2(x, \zeta_2;h) &= e^{\frac{i x \cdot \zeta_2}{h}}(a_2(x,\mu_1 - i \mu_2) + h^{m/2} r_2(x,\zeta_1;h)),
\end{align*}
where the amplitudes $a_1(x,\mu_1 + i \mu_2)$, $a_2(x,\mu_1 - i \mu_2) \in C^{\infty}(\bar B)$ satisfy the transport equations
\begin{equation}\label{actual transport equations 1}
\begin{aligned}
&((\mu_1 + i \mu_2) \cdot \nabla)^{2} a_2(x,\mu_1 + i \mu_2) = 0\quad\mbox{in } B,  \\
&((\mu_1 - i \mu_2) \cdot \nabla)^{2} a_1(x,\mu_1 - i \mu_2) = 0\quad\mbox{in } B,
\end{aligned}
\end{equation}
and the remainder terms $r_1(.,\zeta_1;h)$ and $r_2(.,\zeta_2;h)$ satisfy
$$
||r_j||_{H_{scl}^{m/2}(B)} = \mathcal{O}(1),\quad j=1,2.
$$

We substitute $u_2$ and $v$ in to the (\ref{identity on larger domain}) and get 
\begin{equation}\label{3 terms}
\begin{aligned}
0 &= \frac{1}{h}b^{B}_{\zeta_2 \cdot (A_2 - A_1)}(a_2 + h^{m/2}r_2,e^{i x \cdot \xi}(\bar a_1 + h^{m/2}\bar r_1))  \\
&\qquad + B^{B}_{A_2 - A_1}(a_2 + h^{m/2}r_2,e^{i x \cdot \xi}(\bar a_1 + h^{m/2}\bar r_1)) \\ 
&\qquad + b^{B}_{q_2 - q_1}(a_2 + h^{m/2}r_2,e^{i x \cdot \xi}(\bar a_1 + h^{m/2} \bar r_1)).
\end{aligned}
\end{equation}
Multiply by $h$ throughout and let $h \to 0$ to get
\begin{equation}\label{integral identity for a}
b^{B}_{(\mu_1 + i \mu_2) \cdot (A_2 - A_1)}(a_2,e^{i x \cdot \xi} \bar a_1)=0.
\end{equation}

Let us justify how we get \eqref{integral identity for a}. We use Proposition \ref{estimates for direct problem on Omega} to show
\begin{align*}
&|B^{B}_{A_2 - A_1}(a_2 + h^{m/2}r_2,e^{i x \cdot \xi}(\bar a_1 + h^{m/2}h\bar r_1))| \\ 
&\leq C ||A_1 - A_2||_{W^{-\frac{m}{2}+1,p'}(\mathbb{R}^n)}||a_2 + h^{m/2}r_2||_{H^{m/2}(B)}\,||\bar a_1 + h^{m/2} \bar r_1||_{H^{m/2}(B)} \\
&\leq C (||a_1||_{H^{m/2}(B)} + ||h^{m/2}r_1||_{H^{m/2}(B)})(||\bar a_2||_{H^{m/2}(B)} + || h^{m/2}\bar r_2||_{H^{m/2}(B)}) \\
&\leq C (||a_1||_{H^{m/2}(B)} + ||r_1||_{H^{m/2}_{scl}(B)})\,(||\bar a_2||_{H^{m/2}(B)} + || \bar r_2||_{H^{m/2}_{scl}(B)}) = \mathcal{O}(1).
\end{align*}

Hence
\begin{equation}\label{decay1}
h|B^{B}_{A_2 - A_1}(a_2 + h^{m/2}r_2,e^{i x \cdot \xi}(\bar a_1 + h^{m/2}\bar r_1))| = \mathcal{O}(h).
\end{equation}

We also have for any $m \geq 2$, using Proposition \ref{estimates for direct problem on Omega},
\begin{align*}
&|b^{B}_{q_1 - q_2}(a_1 + h^{m/2}r_1, e^{i x \cdot \xi}(\bar a_2 + h^{m/2} \bar r_2))|  \\
&\leq C ||q_1 - q_2||_{W^{-\frac{m}{2}+\delta,r'}(\mathbb{R}^n)}||a_1 + h^{m/2}r_1||_{H^{\frac{m}{2} - \delta}(B)}\,||\bar a_2 + h^{m/2} \bar r_2||_{H^{\frac{m}{2} - \delta}(B)}  \\
&\leq C (||a_1||_{H^{m/2}(B)} + ||r_1||_{H^{m/2}_{scl}(B)})\,(||\bar a_2||_{H^{m/2}(B)} + ||\bar r_2||_{H^{m/2}_{scl}(B)})  = \mathcal{O}(1).
\end{align*}

Hence,
\begin{equation}\label{decay2}
h|b^{B}_{q_1 - q_2}(a_1 + h^{m/2}r_1, e^{i x \cdot \xi}(\bar a_2 + h^{m/2} \bar r_2))| = \mathcal{O}(h).
\end{equation}
Thus, we see that after multiplying \eqref{3 terms} by $h$, the latter 2 terms in \eqref{3 terms} go to zero as $h \to 0$. 

We also need to justify that
\begin{equation}\label{decay3}
\begin{aligned}
|b^{B}_{\zeta_2 \cdot (A_2-A_1)}(a_2, e^{ix \cdot \xi} &h^{m/2}\bar r_1)|=\mathcal O(h),\quad |b^B_{\zeta_2 \cdot (A_2-A_1)}(h^{m/2}r_2, e^{ix \cdot \xi} h^{m/2}\bar a_1)| = \mathcal O(h),\\
|&b^B_{\zeta_2 \cdot (A_2-A_1)}(h^{m/2}r_2, e^{ix \cdot \xi} h^{m/2}\bar r_1)|=\mathcal O(h)
\end{aligned}
\end{equation}

We only show why $b^{B}_{\zeta_2 \cdot (A_2-A_1)}(a_2, e^{ix \cdot \xi} h^{m/2}\bar r_1)  = \mathcal{O}(h)$. The proof for other two terms follows similarly. By Proposition \ref{MultiplicationTheorem1}, we have for any $m \geq 2$,
\begin{align*}
&|b^{B}_{\zeta_2 \cdot (A_2-A_1)}(a_2, e^{ix \cdot \xi} h^{m/2}\bar r_1)|\\ &\leq ||A_2 - A_1||_{W^{-\frac{m}{2} +1,p'}(\mathbb{R}^n)} ||h^{m/2} \bar r_1||_{H^{\frac{m}{2} -1}(B)} ||a_2||_{H^{\frac{m}{2}}(B_2)} \\
& =\mathcal{O}(h)||r_1||_{H^{\frac{m}{2}-1}_{scl}(B)} = \mathcal{O}(h).
\end{align*}

From \eqref{decay1}, \eqref{decay2} and \eqref{decay3} we see that \eqref{integral identity for a} is indeed justified. 

\begin{remark}\label{rmk::remark33}
	{\it Observe that because our $A$ and $q$ are rough, by duality and Sobolev multiplication, we get estimates in $H^{m/2}(\Omega)$ norm and hence we need a decay of $h^{m/2}$ so that the $H^{m/2}_{scl}$ norm of the correction term is $\mathcal{O}(1)$. If we had just used an $\mathcal{O}(h)$ decay then we would eventually have to use Sobolev estimates in $H^{1}(\Omega)$, which would require $A$ and $q$ to have higher regularity.}
\end{remark}

Now plug in $a_1 = a_2 = 1 $ in (\ref{integral identity for a}) to obtain
$$
\langle (\mu_1 + i \mu_2) \cdot (A_1 - A_2), e^{i x \cdot \xi} \rangle = 0.
$$

We can run the whole argument starting from the construction of $\zeta_1$ and $\zeta_2$, this time with the triple $(\mu_1,-\mu_2,\xi)$, to obtain
$$
\langle (\mu_1 - i \mu_2) \cdot (A_1 - A_2), e^{i x \cdot \xi} \rangle = 0.
$$

The last two equations then imply
\begin{equation}\label{3}
\mu \cdot (\hat A_2(\xi) - \hat A_1(\xi))=0\quad\text{for all}\quad\mu,\xi \in \mathbb{R}^n\quad\text{with}\quad\mu \cdot \xi =0.
\end{equation}
For each $\xi = (\xi_1,\xi_2,...,\xi_n)$ and for $j \neq k$, $1 \leq j,k \leq n$, consider the vector $\mu = \mu(\xi,j,k)$ such that $\mu_j = -\xi_k$, $\mu_k = \xi_j$ and all other components equal to zero.

 Therefore, $\mu$ satisfies $\mu \cdot \xi = 0$. Hence, from (\ref{3}), we obtain
$$
\xi_j \cdot (\hat A_{1,k}(\xi) - \hat A_{2,k}(\xi)) - \xi_k \cdot (\hat A_{1,j}(\xi) - \hat A_{2,j}(\xi)) = 0,
$$

which proves
$$
\partial_{j}(A_{1,k} - A_{2,k}) - \partial_{k}(A_{1,j} - A_{2,j}) = 0 \quad\text{in}\quad \Omega, \quad 1 \leq j,k \leq n,
$$

in the sense of distributions.

To prove $A_1 = A_2$, we consider $A_1 - A_2$ as a 1- current and using the Poincare lemma for currents, we conclude that there is a $g \in \mathcal{D}'(\mathbb{R}^n)$ such that 
$\nabla g = A_1 - A_2$; see \cite{deRham}. Note that $g$ is a constant outside $\bar B$ since $A_1 - A_2 = 0$ in $\mathbb{R}^n \setminus \bar B$ (also near $\partial B$). Considering $g - c$ instead of $g$, we may instead assume $g \in \mathcal{E}'(\bar B)$.

To show $A_1 = A_2$, consider (\ref{integral identity for a}) with $a_1(\cdot,\mu_1 - i \mu_2)=1$ and $a_2(\cdot,\mu_1 + i \mu_2)$ satisfying
$$
((\mu_1 + i \mu_2) \cdot \nabla)a_2(x,\mu_1 + i \mu_2) = 1\quad\text{in}\quad B.
$$

Such a choice of $a_2(\cdot,\mu_1 + i \mu_2)$ is possible because of (\ref{actual transport equations 1}). The previous equation is an inhomogeneous $\bar \partial$-equation and we can solve it by setting
$$
a_2(x,\mu_1 + i \mu_2) = \frac{1}{2 \pi} \int \limits_{\mathbb{R}^2}\, \frac{\chi(x - y_1 \mu_1 - y_2 \mu_2)}{y_1 + i y_2}\,dy_1\,dy_2,
$$

where $\chi \in C^{\infty}_{0}(\mathbb{R}^n)$ is such that $\chi = 1$ near $\bar B$; see \cite[Lemma 4.6]{Mikko}.

From (\ref{integral identity for a}), we have
$$
b^{B}_{(\mu_1 + i \mu_2) \cdot \nabla g}(a_2,e^{i x \cdot \xi})=0.
$$

Now, use the fact that $\mu_1 \cdot \xi = \mu_2 \cdot \xi = 0$ to get 
\begin{align*}
0 &= -b^{B}_{(\mu_1 + i \mu_2) \cdot \nabla g}(a_2,e^{i x \cdot \xi}) = - \langle (\mu_1 + i \mu_2) \cdot \nabla g, e^{i x \cdot \xi}a_2 \rangle_{B}\\
&=\langle g, e^{i x \cdot \xi}(\mu_1 + i \mu_2) \cdot \nabla a_2 \rangle_{B} = \langle g, e^{i x \cdot \xi} \rangle_{B}.
\end{align*}

Since $g$ is compactly supported, this gives $g=0$ in $\mathbb{R}^n$, and in $B$ in particular, implying $A_1 = A_2$.

To show $q_1= q_2$, substitute $A_1 = A_2$ and $a_1 = a_2 = 1$ in to the identity (\ref{identity on larger domain}) to obtain
\begin{align}\label{noproblemtillhere}
b^{B}_{q_2 - q_1}(1 + h^{m/2}r_2,(1 + h^{m/2}\bar r_1) e^{i x \cdot \xi})=0.
\end{align}

Let $h \to 0$ to get $\hat q_1(\xi)- \hat q_2(\xi)=0$ for all $\xi \in \mathbb{R}^n$. To justify this we need to show that
\begin{align*}
b^{B}_{q_2 - q_1}(h^{m/2}r_1,e^{i x \cdot \xi})\to 0,\, b^{B}_{q_2 - q_1}(h^{m/2}r_2,e^{i x \cdot \xi})\to 0,\,\\ b^{B}_{q_2 - q_1}(h^{m/2}r_2,h^{m/2}\bar r_1e^{i x \cdot \xi})\to 0
\end{align*}
as $h\to 0$. We will only consider the term $b^{B}_{q_2 - q_1}(h^{m/2}r_1,e^{i x \cdot \xi})$. The justification for the other two terms follows similarly. We have for any $m \geq 2$ 
\begin{align*}
|b^{B}_{q_2 - q_1}(h^{m/2}r_1,e^{i x \cdot \xi})|  &\leq C ||q_2 - q_1||_{W^{-\frac{m}{2}+\delta,r'}(\mathbb{R}^n)} ||e^{i x \cdot \xi}||_{H^{\frac{m}{2} - \delta}(B)} ||h^{\frac{m}{2}}r_1||_{H^{\frac{m}{2} - \delta}(B)}\\
&=\mathcal{O}(h^{\delta}) ||r_1||_{H^{\frac{m}{2} - \delta}_{scl}(B)}=\mathcal{O}(h^{\delta}) ||r_1||_{H^{\frac{m}{2}}_{scl}(B)}=\mathcal{O}(h^{\delta}).
\end{align*}

Since $\hat q_1(\xi)- \hat q_2(\xi)=0$ for all $\xi \in \mathbb{R}^n$, we get $q_1 = q_2$ in $B$.

\begin{remark}\label{remark::why we need delta positive}
	{\it If we take $\delta = 0$, then we see that all we can say using Proposition~\ref{MultiplicationTheorem1} is that $|b^{B}_{q_2 - q_1}(h^{m/2}r_1,e^{i x \cdot \xi})| = \mathcal{O}(1)$. This is why we impose slightly higher regularity for $q$.}
\end{remark}	

\appendix

\section{Properties of $D_A$ and $m_q$}\label{Appendix}
\vspace{2mm}
The results on the forward problem, as stated in the appendices, hold for any open bounded set $\Omega \subset \mathbb{R}^n$ with smooth boundary, and an integer $m \geq 2$.  

Let $A$ and $q$ satisfy (\ref{defition of A}) and (\ref{defition of q}), respectively. As before, in what follows, $W^{s,p}$ is the standard $L^p$ based Sobolev space on $\mathbb{R}^{n}$, $s \in \mathbb{R}$ and $1 <p < \infty$ defined using Bessel potential.

We start by considering the bi-linear forms 
$$
B_A^{\mathbb{R}^n}(\tilde u,\tilde v) = \langle A, \tilde v D \tilde u \rangle, \quad b_q^{\mathbb{R}^n}(\tilde u, \tilde v) = \langle q , \tilde u \tilde v \rangle, \qquad \tilde u, \tilde v \in H^{m}(\mathbb{R}^n).
$$

The following result shows that the forms $B_{A}^{\mathbb{R}^n}$ and $b_{q}^{\mathbb{R}^n}$ are bounded on $H^{m}(\mathbb{R}^n)$. The proof is based on a property of multiplication of functions in Sobolev spaces.

\begin{proposition} \label{estimates for direct problem on a larger domain}
	The bi-linear forms $B_{A}^{\mathbb{R}^n}$ and $b_{q}^{\mathbb{R}^n}$ on $H^{m}(\mathbb{R}^n)$ are bounded and satisfy for any $m \geq 2$,
	\begin{align*}
	&|b_{q}^{\mathbb{R}^n}(\tilde u, \tilde v)| \leq C ||q||_{W^{-\frac{m}{2}+\delta,r'}(\mathbb{R}^n)}\,||\tilde u||_{H^{\frac{m}{2}}(\mathbb{R}^n)}\,||\tilde v||_{H^{\frac{m}{2}}(\mathbb{R}^n)}, \\
	& |B_A^{\mathbb{R}^n}(\tilde u, \tilde v)| \leq C ||A||_{W^{-\frac{m}{2}+1+,p'}(\mathbb{R}^n)}\,||\tilde u||_{H^{\frac{m}{2}}(\mathbb{R}^n)}\,||\tilde v||_{H^{\frac{m}{2}}(\mathbb{R}^n)}.
	\end{align*}
\end{proposition}
\begin{proof}
	Using the duality between $W^{-\frac{m}{2} + \delta,r'}(\mathbb{R}^n)$ and $W^{\frac{m}{2}-\delta,r}(\mathbb{R}^n)$, we conclude from  Proposition \ref{MultiplicationTheorem1}  that for all $\tilde u, \tilde v \in H^{m}(\mathbb{R}^n)$ with $m \geq 2$,
	\begin{align*}
	|b_{q}^{\mathbb{R}^n}(\tilde u, \tilde v)| & \leq C ||q||_{W^{-\frac{m}{2}+\delta,r'}(\mathbb{R}^{n})}||\tilde u \tilde v||_{W^{\frac{m}{2}-\delta,r}(\mathbb{R}^n)} \\
	&\leq C ||q||_{W^{-\frac{m}{2}+\delta,r'}(\mathbb{R}^{n})} ||\tilde u||_{H^{\frac{m}{2} - \delta}(\mathbb{R}^n)} ||\tilde v||_{H^{\frac{m}{2} - \delta}(\mathbb{R}^n)}\\
	&\leq C ||q||_{W^{-\frac{m}{2}+\delta,r'}(\mathbb{R}^{n})} ||\tilde u||_{H^{\frac{m}{2}}(\mathbb{R}^n)} ||\tilde v||_{H^{\frac{m}{2}}(\mathbb{R}^n)}.
	\end{align*}
	\begin{remark}
	{\it The hypotheses for Proposition \ref{MultiplicationTheorem1} are satisfied for $m \leq n$ with $p_1 = p_2 = 2,\, r=p \in (1, \frac{2n}{2n-m+2 \delta}]$, $s_1 = \frac{m}{2} - \delta, \,s_2 = \frac{m}{2} - \delta$. For $m>n$, we choose $p_1=p_2=2$, $r=p \in (1,2]$, $s_1 = \frac{m}{2} - \delta, \,s_2 = \frac{m}{2} - \delta$.}
	\end{remark}
	
	We now give the estimate for the bi-linear form $B_A^{\mathbb{R}^n}$. Using the duality between $W^{-\frac{m}{2}+1,p'}(\mathbb{R}^n)$ and $W^{-\frac{m}{2}+1,p}(\mathbb{R}^n)$   we conclude from  Proposition \ref{MultiplicationTheorem1} that for all $\tilde u, \tilde v \in H^{m}(\mathbb{R}^n)$, for $m>2$ we have 
	\begin{align*}
	|B_A^{\mathbb{R}^n}(\tilde u, \tilde v)|&\leq C ||A||_{W^{-\frac{m}{2}+1,p'}(\mathbb{R}^n)}\,||D \tilde u \tilde v|_{W^{\frac{m}{2}-1,p}(\mathbb{R}^n)} \\
	&\leq C ||A||_{W^{-\frac{m}{2}+1,p'}(\mathbb{R}^n)}\,||D \tilde u||_{H^{\frac{m}{2} -1}(\mathbb{R}^n)} ||\tilde v||_{H^{\frac{m}{2}}(\mathbb{R}^n)}\\
	&\leq  C ||A||_{W^{-\frac{m}{2}+1,q}(\mathbb{R}^n)}\,||\tilde u||_{H^{\frac{m}{2}}(\mathbb{R}^n)} ||\tilde v||_{H^{\frac{m}{2}}(\mathbb{R}^n)}.
	\end{align*}
	\begin{remark}
	{\it The hypotheses for Proposition \ref{MultiplicationTheorem1} are satisfied for $m<n$  with $p_1 = p_2 = 2$, $p \in (1, \frac{2n}{2n-m}]$, $s_1 = \frac{m}{2} - 1$, $s_2 = \frac{m}{2}$. For $m=n$ and $m=n+2$, we choose $p_1 = p_2 = 2$, $p \in (1, 2)$, $s_1 = \frac{m}{2} - 1$, $s_2 = \frac{m}{2}$. Finally, for other $m$, we choose $p_1 = p_2 = 2$, $p \in (1, 2]$, $s_1 = \frac{m}{2} - 1$, $s_2 = \frac{m}{2}$.}
	\end{remark}
	For the case $m=2$, using H\"older's inequality and Sobolev embedding we get
	\begin{align*}
	|B_A^{\mathbb{R}^n}(\tilde u, \tilde v)| &\leq C ||A||_{L^{n}(\mathbb{R}^n)}\,||D \tilde u \tilde v|_{L^{n'}(\Omega)} \\
	&\leq C ||A||_{L^{n}(\mathbb{R}^n)}\,||D \tilde u||_{H^{\frac{m}{2} -1}(\Omega)} ||\tilde v||_{H^{\frac{m}{2}}(\Omega)} \\
	&\leq  C ||A||_{L^{n}(\mathbb{R}^n)}\,||\tilde u||_{H^{\frac{m}{2}}(\mathbb{R}^n)} ||\tilde v||_{H^{\frac{m}{2}}(\mathbb{R}^n)}.
	\end{align*}
	The proof is thus complete.
\end{proof}

Now we show that the operators $B_A$ and $b_q$ defined in \eqref{bilinearformdefn} are indeed well defined. Recall that
$$
B_A(u,v): = B_{A}^{\mathbb{R}^n}(\tilde u, \tilde v),\quad b_q(u,v): = b_{q}^{\mathbb{R}^n}(\tilde u, \tilde v),\qquad u, v \in H^{m}(\Omega),
$$
where $\tilde u, \tilde v \in H^{m}(\mathbb{R}^n)$ are any extensions of $u$ and $v$, respectively. We want to show that this definition is independent of the choice of extensions $\tilde u, \tilde v$. Indeed, let $u_1, u_2 \in H^m(\mathbb{R}^n)$ be such that $u_1 = u_2 = u$ in $\Omega$, and let $v_1, v_2 \in H^m(\mathbb{R}^n)$ be such that $v_1 = v_2 = v$ in $\Omega$. It is enough to show that for all $w\in H^m(\mathbb R^n)$,
$$
B_A^{\mathbb{R}^n}(u_1,w) = B_A^{\mathbb{R}^n}(u_2, w),  \quad   B_A^{\mathbb{R}^n}(w, v_1) =B_A^{\mathbb{R}^n}(w, v_2)
$$
and 
$$
b_q^{\mathbb{R}^n}(u_1,w) = b_q^{\mathbb{R}^n}(u_2, w),\quad  b_q^{\mathbb{R}^n}(w, v_1) =b_q^{\mathbb{R}^n}(w, v_2).
$$
Since $A$ and $q$ are supported in $\bar \Omega$ and since $u_1=u_2$ and $v_1=v_2$ in $\Omega$, 
we have
$$
B_A^{\mathbb{R}^n}(u_1-u_2,w)=\langle A,D(u_1-u_2)w\rangle=0,\quad b_q^{\mathbb{R}^n}(w,v_1-v_2)=\langle q,w(v_1-v_2)\rangle=0
$$
and
$$
B_A^{\mathbb{R}^n}(w,v_1-v_2)=\langle A,(v_1-v_2)Dw\rangle=0,\quad b_q^{\mathbb{R}^n}(u_1-u_2,w)=\langle q,(u_1-u_2)w\rangle=0.
$$

The next result shows that the bi-linear forms $B_A$ and $b_q$ are bounded on $H^m(\Omega)$.  

\begin{proposition}\label{estimates for direct problem on Omega} The bi-linear forms $B_A$ and $b_q$ are bounded on $H^m(\Omega)$ are bounded and satisfy for any $m \geq 2$
	\begin{align*}
	|b_q(u,v)| &\leq C ||q||_{W^{-\frac{m}{2} + \delta,r'}(\mathbb{R}^n)} ||u||_{H^{\frac{m}{2}}(\Omega)} ||v||_{H^{\frac{m}{2}}(\Omega)} \\
	|B_A(u,v)| &\leq C ||A||_{W^{-\frac{m}{2} + 1 ,p'}(\mathbb{R}^n)}||u||_{H^{\frac{m}{2}}(\Omega)} ||v||_{H^{\frac{m}{2}}(\Omega)}
	\end{align*}
	for all $u,v \in H^m(\Omega)$.
\end{proposition}
\begin{proof}This easily follows from the previous proposition in exactly the same way as in \cite[Proposition~A.2]{Yernat}.
\end{proof}
Now, for $u \in H^m(\Omega)$, we define $D_A(u)$ and $m_q(u)$ for any $v \in H^m_{0}(\Omega)$ by 
$$
\langle D_A(u),v \rangle_\Omega := B_A(u,v), \quad \langle m_q(u), v \rangle_\Omega := b_q(u,v).
$$
The following result, which is an immediate corollary of Proposition {\ref{estimates for direct problem on Omega}, implies that $D_A$ and $m_q$ are bounded operators from $H^m(\Omega) \to H^{-m}(\Omega)$. The norm on $H^{-m}(\Omega)$ is the usual dual norm given by
	\[
	||v||_{H^{-m}(\Omega)} = \underset{0 \neq \phi \in H^m_{0}(\Omega)}{sup} \frac{|\langle v, \bar\phi \rangle_\Omega|}{||\phi||_{H^m(\Omega)}}.
	\]
	
	\begin{corollary}\label{boundedness of A and q}
		The operators $D_A$ and $m_q$ are bounded from $H^m(\Omega) \to H^{-m}(\Omega)$ and satisfy
		\begin{align}
		&||m_q(u)||_{H^{-m}(\Omega)} \leq C ||q||_{W^{-\frac{m}{2}+\delta,r'}(\mathbb{R}^n)}||u||_{H^{m}(\Omega)} \mbox{ and}\notag \\
		&||D_A(u)||_{H^{-m}(\Omega)} \leq C ||A||_{W^{-\frac{m}{2}+1,p'}(\mathbb{R}^n)}||u||_{H^{m}(\Omega)}
		\end{align}
		for all $u \in H^m (\Omega)$.
	\end{corollary}
	Finally, we state the following identities which are useful for defining the adjoint of $\mathcal{L}_{A,q}$.
	\begin{proposition}\label{Adjoint}
		For any $u,v \in H^m(\Omega)$, the forms $B_A$ and $b_q$ satisfy the following identities
		$$
		B_A(u,v) = -B_A(v,u) - b_{D \cdot A}(u,v) \quad and \quad b_q(u,v) = b_q(v,u).
		$$
	\end{proposition}
	\begin{proof}
		Since the proof repeats that of \cite[Proposition~A.4]{Yernat} almost word for word, we omit it.
	\end{proof}
	
	\section{Well-posedness and Dirichlet-to-Neumann map}
	
	Let $\Omega \subset \mathbb{R}^n$, $n \geq 3$ be any bounded open set with smooth boundary, and let $A$ and $q$ be as in (\ref{defition of A}) and (\ref{defition of q}) respectively with $m \geq 2$.
	
	\begin{proposition}\label{well-posedness}
	The Dirichlet-to-Neumann map $\mathcal{N}_{A,q}$ is a bounded operator from $$
	\prod_{j=0}^{m-1} H^{m - j - 1/2}(\partial \Omega) \rightarrow \prod_{j=0}^{m-1} H^{-m + j + 1/2}(\partial \Omega).
	$$
\end{proposition} 
\begin{proof}
The proof follows well-known variational argument principles. 

For $f = (f_0,f_1,....,f_m-1) \in \prod_{j=0}^{m-1} H^{m-j-1/2}(\partial \Omega)$, we consider the Dirichlet problem
	\begin{equation}\label{direct}
	\mathcal L_{A,q}u = 0\text{ in }\Omega\quad\text{and}\quad\gamma u = f \text{ on } \partial \Omega,
	\end{equation}
	where $\gamma $ is the Dirichlet trace operator $\gamma : H^{m}(\Omega) \to \prod_{j=1}
	^{m-1} H^{m-j-1/2}(\partial \Omega)$ which is bounded and surjective; see \cite[Theorem 9.5]{Grubb}.
	
	First aim of this appendix is to use standard variational arguments to show well-posedness of problem (\ref{direct}). We start with the following inhomogeneous problem
	\begin{equation}\label{inhomo}
	\mathcal L_{A,q}u = F\mbox{ in }\Omega\quad\text{and}\quad \gamma u = 0 \mbox{ on }\partial \Omega,\quad F \in H^{-m}(\Omega).
	\end{equation}
	To define a sesqui-linear form $a$ associated to the problem (\ref{inhomo}), for $u,v \in C^{\infty}_{0}(\Omega)$, we can integrate by parts and obtain 
	$$
	\langle\mathcal  L_{A,q}u,\bar v \rangle_{\Omega} = \sum _{|\alpha| = m} \frac{m !}{\alpha !} \int \limits_{\Omega} D^{\alpha} u \,\overline{D^{\alpha} v} \,dx+ \langle D_A(u),  \bar v \rangle_{\Omega} + \langle m_q(u),  \bar v \rangle_{\Omega}:=a(u,v).
	$$
	
	Hence we define $a$ on $H^m_0(\Omega)$ by 
	\begin{align*}
	a(u,v): = \sum _{|\alpha| = m} \frac{m !}{\alpha !}\int \limits_{\Omega} D^{\alpha} u  \,\overline{D^{\alpha} v} \,dx + \langle D_A(u),  \bar v \rangle_{\Omega} + \langle m_q(u), \bar v \rangle_{\Omega},\;\; u,v \in H^m_0(\Omega).
	\end{align*}
	
	We now show that this sesquilinear form $a$ is bounded on $H^m_{0}(\Omega)$. Using duality and  Proposition \ref{estimates for direct problem on Omega}, for $u,v \in H^m_{0}(\Omega)$, we obtain
	$$
	|a(u,v)| \leq C ||u||_{H^m(\Omega)}||v||_{H^m(\Omega)},
	$$
	thereby showing boundedness of $a$. Moreover, Poincare's inequality for $u \in H^m_0(\Omega)$ gives 
	$$
	||u||^2_{H^m(\Omega)} \leq C \sum_{|\alpha| = m} ||D^{\alpha} u||^2_{L^2(\Omega)}.
	$$
	
	Split $q = q^\sharp +(q - q^\sharp)$ with $q^\sharp \in L^{\infty}(\Omega, \mathbb{C})$ and $||q - q^\sharp||_{W^{-\frac{m}{2}+\delta,r'}(\mathbb{R}^n)}$ small enough, and split $A = A^\sharp  + (A - A^\sharp)$ with $A^\sharp \in L^{\infty}(\Omega, \mathbb{C}^n)$ and 
	$||A - A^\sharp||_{W^{-\frac{m}{2}+1,p'}(\mathbb{R}^n)}$ small enough. Using Poincare's inequality and Proposition \ref{estimates for direct problem on Omega}, we obtain,
	\begin{align*}
	&Re\, a(u,u) \\&\geq \sum_{|\alpha|=m } \frac{m !}{\alpha !} ||D^\alpha u||^2_{L^2(\Omega)}  - |B_A(u, \bar u)| - |b_{q}(u,\bar u)| \\
	&\geq C \sum_{|\alpha| = m} ||D^\alpha u||^2_{L^2(\Omega)} - |B_{A^\sharp}(u,\bar u)| - |b_{q^\sharp}(u, \bar u)| - |B_{A - A^\sharp}(u,\bar u)| - |b_{q - q^\sharp}(u,\bar u)| \\
	&\geq C ||u||^2_{H^m(\Omega)} - ||A^\sharp||_{L^{\infty}(\Omega)} ||Du||_{L^2(\Omega)}||u||_{L^2(\Omega)} - ||q^\sharp||_{L^{\infty}(\Omega)}||u||^2_{L^2(\Omega)} \\
	&\quad-C' ||A - A^\sharp||_{W^{-\frac{m}{2}+1,p'}(\mathbb{R}^n)}||u||^2_{H^{\frac{m}{2}}(\Omega)} - C'||q - q^\sharp||_{W^{-\frac{m}{2}+\delta,r'}(\mathbb{R}^n)}||u||^2_{H^{\frac{m}{2}}(\Omega)} \\
	&\geq C ||u||^2_{H^m(\Omega)} -  ||A^\sharp||_{L^{\infty}(\Omega)}\frac{\epsilon}{2}||Du||^2_{L^2(\Omega)} - ||A^\sharp||_{L^{\infty}(\Omega)}\frac{1}{2 \epsilon}||u||^2_{L^2(\Omega)} \\
	&\quad - ||q^\sharp||_{L^{\infty}(\Omega)}||u||^2_{L^2(\Omega)} -C' ||A - A^\sharp||_{W^{-\frac{m}{2}+1,p'}(\mathbb{R}^n)}||u||^2_{H^{\frac{m}{2}}(\Omega)} \\
	&\quad - C'||q - q^\sharp||_{W^{-\frac{m}{2}+\delta,r'}(\mathbb{R}^n)}||u||^2_{H^{\frac{m}{2}}(\Omega)}.
	\end{align*}
	
	Now choose $\epsilon >0$ to be sufficiently small to get
	$$
	Re\,a(u,u) \geq C ||u||^2_{H^m(\Omega)} - C_0||u||^2_{L^2(\Omega)} \quad  C, C_0 >0,\quad u \in H^m_{0}(\Omega).
	$$
	
	Therefore, the seqsuilinear form $a$ is coercive on $H^m_0(\Omega)$. Compactness of the embedding $H^{m}_0(\Omega) \hookrightarrow H^{-m}(\Omega)$ together with positivity of bounded operator $\mathcal{L}_{A,q} + C_{0} I:H^m_{0}(\Omega) \to H^{-m}(\Omega)$ imply that $\mathcal{L}_{A,q}: H^m_{0}(\Omega) \to H^{-m}(\Omega)$ is Fredholm with zero index and hence Fredholm alternative holds for $\mathcal{L}_{A,q}$; see \cite[Theorem~2.33]{Mclean}. (\ref{inhomo}) thus has a unique solution $u \in H^m_0(\Omega)$ if $0$ is outside the spectrum of $\mathcal{L}_{A,q}$.
	
	Now, consider the Dirichlet problem (\ref{direct}) and assume $0$ is not in the spectrum of $\mathcal{L}_{A,q}$. We know that there is a $w \in H^m(\Omega)$ such that $\gamma w = f$. According to the Corollary \eqref{boundedness of A and q}, we have $\mathcal{L}_{A,q}w \in H^{-m}(\Omega)$. Therefore $u = v + w$ with $v \in H^m_{0}(\Omega)$ being the unique solution of the equation $\mathcal{L}_{A,q}v = -\mathcal{L}_{A,q}w \in H^{-m}(\Omega)$ is the unique solution of the Dirichlet problem (\ref{direct}).
	
	Under the assumption that $0$ is not in the spectrum of $\mathcal{L}_{A,q}$, the Dirichlet-to-Neumann map is defined as follows: Let $f,h \in \prod_{j=0}^{m-1} H^{m-j-1/2}(\partial \Omega)$. Set
	\begin{equation}\label{def::Neumann data}
	\langle \mathcal{N}_{A,q}f,\bar h \rangle_{\partial \Omega}: = \sum_{|\alpha|=m}\frac{m !}{\alpha !}   (D^{\alpha }u, D^{\alpha} v_{h})_{L^{2}(\Omega)} + B_A(u, \bar v_h) + b_q(u,  \bar v_h),
	\end{equation}
	where $u$ is the unique solution of the Dirichlet problem (\ref{direct}) and $v_h \in H^m(\Omega)$ is an extension of $h$, that is $\gamma v_h = h$. To see that this definition is independent of $v_h$, let $v_{h,1}, v_{h,2} \in H^m(\Omega)$ be such that $\gamma v_{h,1} = v_{h,2} = h$. Since $w = v_{h,1} - v_{h,2} \in H^m_{0}(\Omega)$ and $u$ solves the Dirichlet problem (\ref{direct}), we have,  
	$$
	0 = \langle \mathcal{L}_{A,q}u,\bar w \rangle_\Omega =\sum_{|\alpha|=m}\frac{m !}{\alpha !}   (D^{\alpha }u, D^{\alpha} w)_{L^{2}(\Omega)} + B_A(u, \bar w) + b_q(u, \bar  w).
	$$
	
	This shows that the definition \eqref{def::Neumann data} is independent of the extension $v_h$. 
	
	Now that we have shown the Dirichlet problem \eqref{direct} is well-posed, we now show that $\mathcal{N}_{A,q}$ is a bounded operator from $$\prod_{j=0}^{m-1} H^{m - j - 1/2}(\partial \Omega) \rightarrow \prod_{j=0}^{m-1} H^{-m + j + 1/2}(\partial \Omega).$$

	From the boundedness of the sesquilinear form $a$ it follows that
	\begin{align*}
	|\langle N_{A,q}f, \bar h \rangle_{\partial \Omega}| &\leq C ||u||_{H^m(\Omega)}||v_h||_{H^m(\Omega)}\\&\leq C ||f||_{\prod_{j=0}^{m-1} H^{m - j - 1/2}(\partial \Omega)} ||h||_{\prod_{j=0}^{m-1} H^{m - j - 1/2}(\partial \Omega)},
	\end{align*}
	where 
	\begin{align*}
	||g||_{\prod_{j=0}^{m-1} H^{-m + j + 1/2}(\partial \Omega)} = (||g_0||^2_{H^{m-1/2}(\partial \Omega)} + ....+||g_{m-1}||^2_{H^{1/2}(\partial \Omega)})^{1/2}
	\end{align*}
	is the product norm on the space $\prod_{j=0}^{m-1} H^{m - j - 1/2}(\partial \Omega)$. Here we have made use of the fact that the extension operator $\prod_{j=0}^{m-1} H^{m - j - 1/2}(\partial \Omega)\ni h \mapsto v_h \in H^m(\Omega)$ is bounded; see \cite[Theorem~9.5]{Grubb}. 
	
	This shows that $\mathcal{N}_{A,q}$ maps $\prod_{j=0}^{m-1} H^{m - j - 1/2}(\partial \Omega)$ into\\ $\big(\prod_{j=0}^{m-1} H^{m - j - 1/2}(\partial \Omega)\big)' =\prod_{j=0}^{m-1} H^{-m + j + 1/2}(\partial \Omega)$ continuously.
\end{proof}	
	
	\section{Bessel potential spaces versus Slobodeckij spaces}
	In this section we show why it is important to consider the Sobolev spaces defined via the Bessel potential. 
	
	There is an alternative, non-equivalent way to generalize the definition of an integer valued Sobolev space to allow fractional exponents. We can define Sobolev spaces with non-integer exponents as Slobodeckij spaces, i.e. if $s = k + \theta$ with $k \in \mathbb{N}_0$ and $\theta \in (0,1)$, then for $p \in [1,\infty)$,
	$$
	H^{s,p}(\mathbb{R}^n) = \{u \in W^{k,p}(\mathbb{R}^n): ||u||_{H^{s,p}(\mathbb{R}^n)}<\infty\},
	$$ 
	where
	\begin{align*}
	||u||_{H^{s,p}(\mathbb{R}^n)} := ||u||_{W^{k,p}(\mathbb{R}^n)} + \bigg(\sum_{|\alpha|=k} \int \int_{\mathbb{R}^n \times \mathbb{R}^n} \frac{|\partial ^{\alpha}(x) - \partial^{\alpha}(y)|^p}{|x-y|^{n+\theta p}} \,dx\,dy\bigg)^{1/p}.
	\end{align*}
	Slobodeckij spaces are special cases of Besov spaces, see \cite{Triebel}. If $s<0$ and $p \in(1,\infty)$, we define 
	$H^{s,p}(\mathbb{R}^n) = (H^{-s,p/(p-1)}(\mathbb{R}^n))^{*}$.
	
	
	
	
	We use the Bessel potential definition in this paper as that definition gives more flexibility with regards to  multiplication as the following result shows.
	\begin{proposition}
		Suppose $s, s_{1} \geq 0$, $s \notin \mathbb{Z}$ and $p_1,p_2, p>1$. If the pointwise multiplication of functions is a continuous bi-linear map $H^{s,p_1}(\mathbb{R}^n) \times H^{s_1,p_1}(\mathbb{R}^n) \hookrightarrow H^{s,p}(\mathbb{R}^n)$, then $p_1 \leq p$.
	\end{proposition}
	\begin{proof}
		Follows from \cite[Proposition~4.3]{Multiplicationpaper}.
	\end{proof}

\end{document}